\documentclass[a4paper,11pt]{amsart}
\usepackage{amscd}
\usepackage{amssymb}
\usepackage{amsthm}
\usepackage{epsf}
\newtheorem{thm}{Theorem}[section]
\newtheorem{cor}[thm]{Corollary}
\newtheorem{exa}[thm]{Example}
\newtheorem{lemma}[thm]{Lemma}
\newtheorem{prop}[thm]{Proposition}
\newtheorem{example}[thm]{Example}
\newtheorem{defn}[thm]{Definition}

\newtheorem{rem}[thm]{Remark}
\newtheorem{conj}[thm]{Conjecture}

\def\min{\operatorname{min}}

\def\max{\operatorname{max}}

\def\c1{\operatorname{c_1}}
\def\c2{\operatorname{c_2}}

\def\ZZ{{\mathbf Z}}

\def\PP{{\mathbf P}}

\def\L{{\mathcal L}}

\def\N{{\mathcal N}}

\def\P{{\mathcal P}}
\def\G{{\mathcal G}}

\def\+{\oplus}                   % direct sum
\def\*{\otimes}                  % tensor product
\def\khpil{\rightarrow}

\begin{document}

\title{Schubert unions in Grassmann varieties}
\author{Johan P. Hansen, Trygve Johnsen and Kristian Ranestad}  

\address{Johan P. Hansen\\ Dept. of Mathematics\\ University of Aarhus, 
Bygn. 530\\ DK-8000 C Aarhus, Denmark}
\email{matjp@imf.au.dk}
\address{Trygve Johnsen\\ Dept. of Mathematics\\ University of Bergen, Johs. Bruns gt 12,
\\ N-5008 Bergen, Norway}
\email{johnsen@mi.uib.no}
\address{Kristian Ranestad\\ Dept. of Mathematics\\ University of Oslo, P.O. 1053,\\N-316 Oslo, Norway}
\email{ranestad@math.uio.no}

%\address{Dept. of Mathematics\\ 
% University of Bergen\\ Johs. Brunsgt 12\\ 5008 Bergen, Norway}
%\email{johnsen@mi.uib.no, andreask@mi.uib.no}
\begin{abstract}
We study subsets of Grassmann varieties $G(l,m)$ over a field $F$, such that these
subsets are unions of Schubert cycles, with respect to a fixed flag.
We study the linear spans of, and in case of positive characteristic, 
the number of $F_q$-rational points on such unions. Moreover we study
a geometric duality of such unions, and give a combinatorial interpretation 
of this duality.
We discuss the maximum number of $F_q$-rational points for Schubert
unions of a given spanning dimension, and we 
give some applications to coding theory. We define Schubert union
codes, and study the parameters and support weights of these codes 
and 
of the well-known Grassmann codes.
\end{abstract}
\maketitle
%\centerline{Addresses: } 
%\centerline{Johan P. Hansen,} 
%\centerline{Dept. of Mathematics,} 
%\centerline{University of Aarhus, Bygn. 530}
%\centerline{DK-8000 C Aarhus, Denmark}
%\centerline{and}
%\centerline{Trygve Johnsen,}
%\centerline{Dept. of Mathematics,} 
%\centerline{University of Bergen, Johs. Bruns gt 12,}
%\centerline{N-5008 Bergen, Norway}
%\centerline{and}
%\centerline{Kristian Ranestad,} 
%\centerline{Dept. of Mathematics,} 
%\centerline{University of Oslo, P.O. 1053,}
%\centerline{N-316 Oslo, Norway}
%\medskip
%\centerline{E-mail: matjp@imf.au.dk, johnsen@mi.uib.no, 
%ranestad@math.uio.no}
%\bigskip
%  \begin{abstract}
%       We give a simple result about some higher weights of Grassmann codes.
%We give some upper bounds for other higher weights, using Schubert cycles.
%We show that there are Schubert cycles that do not compute higher weights.
%At the ende we consider some natural questions.
%  \end{abstract}
%
%\vspace{.75cm}

% Section 1
\section{Introduction}
\label{intro}

Let $G(l,m)=G_F(l,m)$ be the Grassmann variety of $l$-dimensional subspaces of a fixed $m$-dimensional vector space $V$ over a field $F$. 
 By the standard Pl\"ucker 
coordinates $G(l,m)$ is embedded into
$\PP^{k-1}=\PP_F^{k-1}$ as a non-degenerate smooth subvariety, where $k=$$m
\choose\ l$. 

This paper was motivated by the problem: ''What is the maximal 
intersection of a linear subspace of a given dimension in $\PP^{k-1}$ 
with $G(l,m)$?'' Here ``maximal'' may refer to the number of $F_{q}$-rational points, 
or to Krull dimension.  We answer the last question, concerning Krull dimension, in Theorem 
\ref{maxkrull}, while for the former we only give partial answers.
The unions of Schubert cycles with respect to a fixed flag, called Schubert unions, turn out 
to play a key role in this problem, so the main body of this paper 
concerns them with the hope that this treatment may have some independant 
interest.

First we point out that the spanning dimension and number of $F_q$-rational points of the Schubert unions are natural generalizations of the corresponding formulas for
Schubert cycles, and that each union is precisely the intersection of its linear span with $G(l,m)$. 

Furthermore we identify the Schubert unions as fixed points for the action of a certain Borel group acting on natural (derived) Grassmann varieties.

This is instrumental in applying Borel's Fixed Point Theorem to solve 
the maximality problem for the Krull dimension.

The set of Schubert unions enjoy a involution given by a natural duality.
The dual of a Schubert union $U$ arises from taking the set of 
$(m-l)$-dimensional subspaces of $V$ intersecting all $l$-dimensional subspaces parametrized by points of $U$ non-trivially. By the standard duality between $G(l,m)$ and $G(m-l,m)$ we interpret the
dual set as a subset of $G(l,m)$. It then turns out that the dual of a Schubert union is a Schubert union with complementary spanning dimension in the Pl\"ucker space (and that we have biduality justifying the terminology).

We proceed to study a natural point grid in ${\bf Z}$$^l$ corresponding to
the Pl\"ucker coordinates. It turns out that the duality just described corresponds to a duality of subsets of this grid, where the cardinality of a subset 
$G_U$ corresponding to $U$ is equal to the spanning dimension of the affine cone of $U$ in the affine space (of dimension $m \choose k$) over the Pl\"ucker space. In the special case $l=2$ we describe an additional structure on the set of
Schubert unions for fixed $m$, as a power set $\P(M)$ for a set $M$ 
with $m-1$ elements. 

In the case where $F$ is the finite
field $F_q$, we study the (equivalence class of) Grassmann code(s) $C(l,m)$ 
obtained by using the Pl\"ucker coordinate tuples of each point of $G(l,m)$ as columns of the generator matrix.
The weight hierarchy of these codes is in a natural manner computed by 
linear sections of $G
(l,m)$, and finding the weights amounts to finding the
linear sections, for each spanning dimension, with the maximum number of $F_q$-rational points. 
It then easily turns out that the highest weights, corresponding to sections of small spanning dimension, are computed by projective spaces that are Schubert cycles and therefore Schubert unions. Furthermore the lowest weights, found in \cite{N}, corresponding to linear sections with large spanning dimension, are computed by the respective dual Schubert unions.
This observation makes it natural to ask whether the linear sections, for each spanning dimension, with the maximum number of $F_q$-rational points, can be taken to be Schubert unions. Moreover we ask whether the dual $U^*$ of a Schubert union $U$ with the maximum number of $F_q$-rational points (for a given spanning dimension) has the maximal number of $F_q$-rational points for the complementary spanning dimension.
It turns out that these two statements are true for certain low values of $l$ and $m$. On the other hand we show that for $l=2$, both statements cannot hold simultaneously, for large enough $m$. In fact we show that if $m$ is large enough, then 
at least one of them fails for spanning dimensions between ca. 36\% and ca. 64\%
of $m \choose\ 2$. In contrast, the percentage of number of weights described in \cite{N}, will approach zero.  For these code-theoretical
questions associated to $G(l,m)$, see also \cite{N}, \cite{GL}, and 
\cite{GT}.

%but in the following sections 
%\ref{SUdef}-\ref{twodim}
In general $F$ will be just a field, unless otherwise specified. The paper is organized as follows.

In Section \ref{SUdef} we will fix an ordered basis for an
$m$-dimensional vector space $V$, and recall the well-known 
cell-decomposition of $G(l,m)$ and the definition of Schubert cycles
with respect to this flag, as described by many other authors.
Then we will proceed to study unions of such cycles, and determine the
dimensions of their linear spans in the Pl\"ucker space  $\PP^{k-1}$, 
and
find concrete equations for their linear spans.
For the case $F=F_q$ we will use the cell decomposition to
determine the number of $F$-rational points on the given unions.

In Section \ref{dual} we will describe the geometric duality of Schubert
unions, and show how this duality corresponds to a quite natural
duality or symmetry of a certain ``diagram of $l$-cubes (or just 
points)''. 
Each $l$-cube in the diagram corresponds to a Pl\"ucker coordinate 
$X_{i_1,..,i_l}$, and the center of the $l$-cube is located  in the 
point $(i_1,..,i_l)$. In the case $F=F_q$ the number of $F$-rational
points on a given union will be easily obtainable from that of its 
dual
union. Our duality, which interchanges dimensions of linear spans of, 
and not Krull dimension of, Schubert unions, is different from, but 
closely related to, Poincare duality.

In Section  \ref{twodim} we will study the particular case $l=2$. We 
will
show how the set of Schubert unions for a fixed $m$ in a natural way
corresponds to a Boolean algebra $P(M)$, where $M$ is a set of $m-1$
elements, say $\{1,2,...,m-1\}$. Moreover duality of Schubert unions
will correspond to complementarity of subsets of $M$, and there will
be an easy way to describe the ``diagram of $2$-cubes''(squares) from 
a subset
of $m$ and vice versa. 

In Section  \ref{codes} and \ref{w67} we will discuss applications to 
coding theory.
The Schubert unions will be used to tell us as much 
as
possible about the support weights of the Grassmann codes $G(l,m)$  
described for example in  \cite{N}, \cite{GL}, and \cite{GT}.
We introduce the upper Schubert union bound for the support weights of
these codes, and show how this bound gives the true value
in two ranges, that described by Nogin in  \cite{N}, and that 
described by Schubert unions that are dual to those used to reproduce
Nogin's result.
We will discuss some natural questions for the 
properties of the support weights and the properties of the Schubert
union bounds, and investigate these properties for the $C(2,m)$
for some low values of $m$. In particular we will determine all ten 
support weights for $C(2,5)$ and show how they coincide with the
Schubert union bound(s).

In Section \ref{algo} we assume $F=F_q$ and $l=2$. We describe a 
systematic  
algorithm for finding the maximum number of $F_q$-rational points on 
a Schubert union spanning an $K$-space in the Pl\"ucker space, at 
least for
large $q$. This amounts 
to finding the Schubert union bound for $d_r$ for $r=k-K$ for each 
$K$ in the range $0 \le K \le k$, where $ k= $$m \choose\ 2$ as 
usual. 
It turns out that for $m \ge 10$ the optimal codes for the different
spanning dimensions do not always exhibit the same symmetries as for 
$m \le
9$. We reveal the asymptotic behaviour of the optimal Schubert unions
of $G(2,m)$ as $m$ goes to infinity.

In Section \ref{Krull} we treat the question: ``What is the maximal 
Krull dimension of a component of an intersection of $G(l,m)$ with a 
linear space of dimension $K$ in the Pl\"ucker space, for each 
possible $K$ ?'' 
%This is a question related to the issue of finding 
%linear intersections (over $F_q$) with the maximal number of 
%$F_q$-rational points, which again is essential for finding higher 
%weights for the Grassmann codes. 
We use Borel's Fixed Point Theorem, and show that the 
maximal Krull dimension is always attained by linear intersections 
spanned by Schubert unions. So far we have not been able to
prove the corresponding fact for the maximal number of $F_q$-rational 
points, (except for the smallest and the largest linear spaces).

In Section \ref{unionskoder} we treat codes made from the points of a
Schubert union. We determine the minimum distance and some higher 
weights of  
such codes, in the case $l=2$. Our results rely heavily on a result 
from \cite{HC}, where the minimum distances for codes from Schubert
cycles in $G(2,m)$ are found.

In the appendix of Section \ref{Tables} we list some known facts 
abouts $d_r$ 
for 
the $G(2,m)$ for some low $m$, and demonstrate how the Schubert unions
of $G(2,m)$ for $m\le 9$, and also $G(3,6)$, exhibit some nice 
symmetries.

We thank Rita Vincenti for an enlightening correspondence at the 
start of this work,  introducing us to \cite{MV}, and thereby 
enabling us to prove Proposition \ref{d5}.
We also thank Torsten Ekedahl for suggesting to use Borel's fixed 
point theorem to prove the results in Section \ref{Krull}.
The second author was supported in part by the Norwegian Research 
Council 
and thanks for this support and for the kind hospitality extended by 
Aarhus University, Denmark, where the support was spent.

\section{Basic Description of Schubert Unions}
\label{SUdef}

In this section we will recall the well known definition of a Schubert
cycle $\alpha =(a_1,...,a_l)$ in the Grassmann variety $G(l,m)$ over 
a field
$F$, and describe unions of such cycles. 
Let $B=\{e_1,...,e_m\}$ be a basis of a $m$-dimensional 
vector space $V$ over $F$, and let $A_i=Span$$\{e_1,...,e_i\}$ in 
$V$, 
for $i=1,...,m.$  Then $A_{1}\subset A_{2}\subset \ldots \subset 
A_{m}=V$ form a complete flag of subspaces of $V$.   With respect to 
the basis $B$ there 
is the following canonical cell decomposition of $G(l,m)$.

For a given $l$-subspace $W$ of $V$
form an $(l \times m)$-matrix $M_{W}$ where the rows form a set of 
basis
vectors for $W$, each row expressed in terms of the basis $B$. 

We choose a basis 
%$<w_{1},\ldots,w_{l}>$ 
for $W$ such that the 
matrix $M_{W}$  have reduced lower left triangular form, i.e. 
the last nonzero entry in each row is $1$, each of these $1$'s are 
the 
only nonzero entries in their column, and each of these $1$'s lie in 
a column to the right of the trailing $1$ in the previous row. 
The trailing $1$ in row $i$ is then in column $a_{i}(W)$ where 
$$a_{i}(W)={\rm min}\{j|\dim(W)\cap A_{j}=i\}.$$
Obviously $$1\leq a_{1}(W)<a_{2}(W)<...<a_{l}(W)\leq m.$$
For $\alpha=(a_1,....,a_l)$ with $ 1 \le a_1 < a_2 <....a_l \le    m$ 
let 
$$C_{\alpha}=\{[W]\in G(l,m)| a_{i}(W)=a_{i}, i=1,...,l\}$$
Since the reduced lower left triangular form of $M_{W}$ is unique,  
$C_{\alpha}$ is an affine space of dimension 
$$\dim C_{\alpha}=\sum_{i=1}^l (a_{i}-i)=\sum_{i=1}^l 
a_{i}-\frac{l(l+1)}{2}.$$
Therefore $C_{\alpha}$ is called a cell, and when $\alpha$ varies 
these cells 
are pairwise disjoint and form a decomposition of $G(l,m)$.

Note that the ordered $l$-uples $\alpha$ belong to the grid 
$$G_{G(l,m)}=\{\beta=(b_1,....,b_l)\in {\bf Z}^l | 1 \le b_1 < b_2 
<....b_l \le    m\}$$
and that this grid is partially ordered by
$\alpha\leq \beta$ if $a_{i}\leq b_{i}$ for $i=1,...,l$.

For each $\alpha\in G_{G(l,m)} $ the Schubert cycle $S_{\alpha}$ is 
defined as:

$$S_{\alpha}=\{W |\dim(W \cap A_{a_i}) \ge i, \quad 
i=1,...,l\}=\cup_{\beta\leq\alpha}C_{\beta}.$$
Note that $S_{\alpha}$ inherits a cell-decomposition from $G(l,m)$.

Next, we choose coordinates for the Pl\"ucker space 
$\PP^{k-1}=\PP(\wedge^l V)$, with
respect to the chosen basis $B$. 
Our choice of Pl\"ucker coordinates are the maximal minors 
of the matrix $M_{W}$ (with alternating signs). The Pl\"ucker 
coordinates from another basis for $W$ 
would differ from these only by a nonzero scalar factor, so they are 
welldefined as projective coordinates.  
Since the maximal minors of $M_{W}$ are indexed by the grid 
$G_{G(l,m)}$ these Pl\"ucker coordinates are 
denoted by $\{X_{\alpha}(W)| \alpha\in G_{G(l,m)}\}$.
 For the unique basis of $W$ above we get $X_{\alpha}=1$ when $[W]\in 
C_{\alpha}$.

These Pl\"ucker coordinates for subspaces belonging to the Schubert 
cycle $S_{\alpha}$ become particularly simple:

\centerline{$S_{\alpha}=\{p \in G(l,m) | X_{\beta}=0$ for all
  $\beta$ with $b_i > a_i$ for some $i \}.$}
We therefore collect all the indices of nonzero Pl\"ucker coordinates 
on the Schubert cycle $S=S_{\alpha}$ in a subgrid of $G_{G(l,m)}$, 
called the
$S$-grid or $\alpha$-grid:
 
\begin{defn}
 $G_S=G_{\alpha}=\{\beta\in G_{G(l,m)}|\beta\leq \alpha\}.$
\end{defn}
In the special case  $\alpha = (m-l+1,...,m-1,m)$, we get
$S_\alpha=G(l,m)$ and $G_{\alpha}=G_{G(l,m)}$.

 \begin{defn} A subset $M\subset G_{G(l,m)}$ is called {\it Borel 
fixed} if it 
 enjoys the property that $\beta\in M$ 
 whenever $\alpha\in M$ and $\beta\leq\alpha$. 
 \end{defn}
 Clearly $G_{\alpha}$ is Borel fixed for each 
 Schubert cycle $S_{\alpha}$.  Notice that $M$ is Borel fixed if and 
only if the union 
 of the cells $\cup_{\alpha\in M}C_{\alpha}$ is closed.
 
For each $\alpha$ we also define:
   $$H_{\alpha}=G_{G(l,m)} - G_{\alpha},$$
in other words the complement of the $\alpha$-grid in the entire
$G(l,m)$-grid.
We then see that 

$$S_{\alpha}=\{[W] \in G(l,m) | X_{\beta}(W)=0, \forall \beta\in 
H_{\alpha}\},$$ 
\begin{itemize} \label{linspan}

\item For a subset $M$ of $G(l,m)\subset \PP (\wedge^l V)$, let 
$\L(M)$ be its linear span in
the projective  Pl\"ucker space $\PP (\wedge^l V).$

\item Moreover we denote by $L(M)$ the linear span of the
affine cone over $M$ in the affine cone over the Pl\"ucker space.

\item For a subset $N$ of an affine space we denote by $L(N)$ its 
linear span. 
\end{itemize}
%For a Schubert cycle $S_{alpha}$ (it is well known) that
%$\dim L(S_{\alpha})$ is equal to
%  the cardinality 
%of $G_{\alpha}$. See for
%example \cite{GL} (In projective Pl\"ucker space the dimension of the
%span is of course one less, but for notational reasons it seems more
%  convenient to operate in the underlying affine space).
  
  The following is well known (see for example \cite{GT}):
  \begin{prop} \label{Sp}  Let $\alpha\in G_{G(l,m)}$ 
  and consider the Schubert cycle $S_{\alpha}\subset G(l,m)\subset 
\PP(\wedge^l V)$ in its 
  Pl\"ucker embedding.
  \begin{enumerate}\item $\L (S_{\alpha})\cap G(l,m)=S_{\alpha}$
  \item $\dim L(S_{\alpha})$ is equal to
  the cardinality 
of the grid $G_{\alpha}$
\item The number of $F_q$-rational points on $S_{\alpha}$ is 
 $\Sigma_{(x_1,...x_l) \in G_{\alpha}} q^{x_1+...+x_l-l(l+1)/2}$.
 \end{enumerate}
\end{prop}
\begin{proof}  The first statement is already proven above.
  The span of $S_{\alpha}$ 
 has codimension equal to  the cardinality of $H_{\alpha}$, so the 
affine dimension of the span equals the cardinality of the complement
 $G_{\alpha}$ in $G_{G(l,m)}$.  The final statement follows 
immediately from the cell-decomposition of $S_{\alpha}$
 (cf. also Theorem 1 and Proposition 11 of \cite{GT}). 
 \end{proof}
 
We will consider finite intersections and finite unions of such 
Schubert cycles
$S_{\alpha}$ with respect to our fixed flag.
Set $\alpha_i=(a_{(i,1)},a_{(i,2)},....,a_{(i,l)})$, for $i=1....,s$. 
From the definition above it is then clear that:
  $\cap_{i=1}^s S_{\alpha_i} = S_{\gamma}$,
where $\gamma=(g_1,....,g_l)$, and $g_j$ is the minimum of the
set $\{a_{1,j},a_{2,j},...,a_{s,j}\},$ for $j=1,....,l.$  Thus the 
intersection of a 
finite set of Schubert cycles $S_{\alpha}$ is again a Schubert cycle.
In particular $\dim L(\cap S_{\alpha_i})$ is equal to
the cardinality of  $G_{\gamma}$.

For a union $S_U=\cup_{i=1}^s S_{\alpha_i}$ of Schubert cycles, 
denote by 
$G_U$ the union $G_U=\cup_{i=1}^sG_{\alpha_i}$, and set 
$H_U=G_{G(l,m)}-G_U$.
Since all Schubert cycles $S_{\alpha}$ has a decomposition of cells 
$C_{\beta}$,
 and all finite intersections of these Schubert cycles are  again 
Schubert cycles, the union $S_{U}$ 
 also has a cell-decomposition inherited from $G(l,d)$:
$$S_{U}=\cup_{\alpha\in G_{U}}C_{\alpha}$$

\begin{prop} \label{int-un}
Let $S_{\alpha_1}, ....,S_{\alpha_s}$ be 
finitely many Schubert cycles with respect to our fixed flag.  Let 
$S_{\gamma}=\cap_{i=1}^s S_{\alpha_i}$  
be their intersection, and let $S_U=\cup_{i=1}^s S_{\alpha_i}$ be 
their union.
\begin{enumerate}
\item The intersection   $S_{\gamma}$
is itself a Schubert cycle with $S$-grid $G_{\gamma}=\cap_{i=1}^s 
G_{\alpha_i} $.
\item  $\L(S_U)\cap G(l,m)=S_{U}.$
\item dim $L(S_{U})$ equals the cardinality of the grid $G_{U}$.
\item  The number of $F_q$-rational points on $S_{U}$ is 
 $\Sigma_{(x_1,...x_l) \in G_{U}} q^{x_1+...+x_l-l(l+1)/2}$.
 \end{enumerate}
 \end{prop}
 
 \begin{proof}   The first statement is already proven above. 
 The ideal  of the linear span of the union $S_{U}$ is defined by the 
intersection 
 of the ideals of the $S_{\alpha_{i}}$.   Therefore 
$\{X_{\beta}|\beta\in H_{U}\}$ 
 form the linear generators of the ideal $I_{U}$ of $S_U$.   
Looking more closely at $H_{U}$, notice that if $\alpha\in H_{U}$, 
then $\beta\in H_{U}$ 
whenever $\alpha\leq\beta$.   We interpret this condition on 
the Pl\"ucker coordinates of a $l$-dimensional subspace $W\subset V$, 
coming from 
a reduced lower left triangular form for the matrix $M_{W}$ as above.
We need the following characterisation of the coordinates 
$X_{\alpha}$ for $\alpha\in H_{U}$:
\begin{lemma} Let $W\subset V$  and let 
$$\{(w_{ij})| 1\le i\le l, 1\le j\le m\}$$ be the entries of the 
matrix $M_{W}$ in reduced lower triangular form. 
Then $X_{\alpha}(W)=0$ for each $\alpha=(a_{1},...,a_{l})\in H_{U}$ 
if and only if   
one of the 
following is satisfied
\begin{enumerate}
    \item $w_{1j}=0$ when $j\ge a_{1}$
    \item $w_{1j}=w_{2j}=0$ when $j\ge a_{2}$
    
    \ldots
    
    \item $w_{ij}=0$ for $i=1,...,l$ when $j\ge a_{l}$
    \end{enumerate}
    \end{lemma}
    \begin{proof} None of the criteria are satisfied, if and only if 
the 
    trailing $1$'s in each row appear in columns 
    $(b_{1},\ldots,b_{l})$ where $b_{i}\ge a_{i}$ for each $i$. But 
    equivalently the Pl\"ucker coodinate $X_{\beta}(W)=1$ for 
    $\beta=(b_{1},\ldots,b_{l})\in H_{U}$, contrary to our 
    assumption.\end{proof}
Now, each itemized condition in the lemma is the condition that
$[W]$ belong to a Schubert cycle: 

\centerline{$w_{1j}=\ldots=w_{rj}=0$ implies that 
dim$W\cap A_{j-1}\ge r$.}

Therefore the collection of conditions 
given by  the element $\alpha=(a_{1},...,a_{l})\in H_{U}$ implies 
that $[W]$ belongs to a union of Schubert cycles with respect to the 
given flag.   Since 
the intersecton of two Schubert cycles with respect to the given flag 
is 
a Schubert cycle, the intersection of a union of Schubert cycles is 
again a Schubert union.  
In particular, it is now straight forward to check that if 
$X_{\beta}(W)=0$ for each $\beta\in H_{U}$, then $[W]$ belongs to the 
Schubert union $S_{U}$, i.e. the linear span of $S_{U}$ intersects 
$G(l.m)$ precisely in $S_{U}$.

The cardinality of $H_{U}$ is the codimension of  $\L(S_{U})$, so the 
dimension of $L(S_{U})$ 
 equals the cardinality of the complement $G_{U}$.
 Finally the number of  $F_{q}$-rational points is counted using the 
cell-decomposition.  Note that one may also use the 
 exclusion-inclusion principle using Proposition \ref{Sp}, 
 since all intersections are Schubert cycles.
 \end{proof}
 
 Notice that when $S_{U}$ is a Schubert union, then the grid $G_{U}$ 
is 
 Borel fixed.  In fact
 \begin{prop} \label{Borel1} The Borel fixed subsets of the grid 
$G_{G(l,m)}$ are 
 precisely the grids $G_{U}$ of Schubert unions.  Similarly, the 
closed unions of cells $C_{\alpha}$ are precisely the Schubert unions.
     \end{prop}
     \begin{proof}  Let $M$ be a Borel fixed  subset of 
$G_{G(l,m)}$.  
    Since $M$ is finite, it has finitely many maximal elements, 
    $\alpha_{1},\alpha_{2},\ldots, \alpha_{r}$ say.  Then 
    $M=\cup_{i=1}^r G_{\alpha_{i}}$.  In particular $M$ is the grid 
of 
    a Schubert union.  Similarly, the second statement  follows  
considering the fact that
   the cell  $C_{\beta}$ lies in the closure of the cell 
$C_{\alpha}$, precisely when $\beta\leq \alpha$.\end{proof}

    % \begin{rem}\label{Borel}
    The name ``Borel fixed'' originates from the a natural action of a
    Borel subgroup $B$ on the Pl\"ucker space.  Extend the scalars
    of $V$ to the algebraic closure $\overline{F_{q}}$ of the field
    $F_{q}$.  Then  $GL(m, \overline{F_{q}})$ acts (from the left) on 
$V$,
on the exterior
    product $\wedge^l V$ and on
    $G(l,m)$.  Let $B\subset GL(m, \overline{F_{q}})$ be the
    subgroup of upper triangular matrices with respect to the basis
    $e_{1},\ldots, e_{m}$ of $V$. Then $B$ is precisely the subgroup
    that fixes the given flag $A_{1}\subset A_{2}\subset\ldots\subset
    A_{m}$.  Let $X=\{X_{\alpha}\}$ be a set of Pl\"ucker coordinates,
    and let $M(X)=\{\beta|X_{\beta}\in X\}$ be the corresponding grid
    in $G_{G(l,m)}$.  Then the linear span of $X$ is fixed by $B$ if
    and only $M(X)$ is Borel fixed:
 %  \end{rem}
    Let $V_{l}=\wedge^l V$.  Then $V_{l}$ is a $k={m\choose
    l}$-dimensional vector space.  Let $r\leq k$ and consider
    $V_{l,r}=\wedge^rV_{l}$, and the Grassmannian $G(r,V_{l})$.
The linear action of $B$ on $V$ and $V_{l}$, clearly induces a
    linear action on $V_{l,r}$ and on $G(r,V_{l})$.  Let
    $$\{e_{\alpha}=e_{a_{1}}\wedge e_{a_{2}}\wedge \ldots \wedge
    e_{a_{l}} |1\leq a_{1}<a_{2}\ldots \leq a_{l}\}$$ be the basis of
    $V_{l}$ with Pl\"ucker coordinates $X_{\alpha}$.  Order this
    basis lexicographically.  Then
    $$\{e_{\alpha_{1}}\wedge e_{\alpha_{2}}\wedge \ldots
    \wedge e_{\alpha_{r}}|(1,2,\ldots,l)\leq
    \alpha_{1}<\alpha_{2}<\ldots <\alpha_{r}\leq (m-l+1,\ldots,m)\}$$ 
form
a basis
    for $V_{l,r}$.
    The only subspaces of $V$
    that are stable under $B$ are the subspaces $A_{i}$, i.e. those
    spanned by $e_{1},e_{2},\ldots,e_{i}$ for some $i$. Therefore the
    only subspaces of $V_{l}$ that are stable under $B$ are those
    spanned by  $\{e_{\alpha}|\alpha\in I\}$, for some finite index
    set $I\subset G_{G(l,m)}$, with the property that
    $\beta=(b_{1},\ldots,b_{l})\in I$ whenever
    $\alpha=(a_{1},\ldots,a_{l})\in I$ and $b_{i}\leq a_{i}$ for all
 $i$. But these index sets are precisely the Borel fixed subsets
    of $G_{G(l,m)}$.  On the other hand, the stable subspaces of
    $V_{l}$ of dimension $r$ are the only fixed points on $G(r,V_{l})$
    under the action of $B$.

\begin{prop} \label{Borel2}  
The Schubert unions of spanning dimension $r$
    define the fixed points under the action of the Borel subgroup
    $B\subset GL(m,F)$ on $G(r,V_{l})$.
\end{prop}

%   \end{rem}
 
\subsection{Duality of Schubert Unions}
\label{dual}

There are various natural forms of duality for
Schubert unions.  One is combinatorial and one is geometric.  
We will define both and show that they coincide.

 First we will describe a geometrical
duality, valid for Schubert unions, but not for general linear
sections of the Grassmann variety $G(l,m)$.

This duality is a restriction of ordinary duality between the 
Pl\"ucker 
space $\PP(\wedge^l V)$ and its natural dual $\PP(\wedge^l V^{*})$.  
Just as $G(l,m)$ has a natural embedding by Pl\"ucker coordinates in 
$\PP(\wedge^l V)$, the Grassmannian $G(l,m)^*$ 
parametrizing rank $l$-subspaces of $V^*$, i.e. $l$-dimensional 
subspaces 
of linear forms on $V$, has an embedding in $\PP(\wedge^l V^*)$.

For a subspace $L\subset  V$, we set 
$$L^{\bot}=\{ H\in V^{*}|L\subset {\rm ker}H\}.$$

Let $M$ be a subset of $G(l,m)$. Then we define:

\begin{defn} \label{duality}
The geometric dual set $D(M)$ is the subset of $G(l,m)^*\subset 
\PP(\wedge^l V^*)=\check{\PP}(\wedge^l V)$
parametrizing $l$-subspaces of linear forms on $V$ whose common 
kernel simultaneously intersect 
all $l$-spaces respresented by points of $M$ in more than $0$.
\end{defn}  

In terms of natural duality between $\PP(\wedge^l V)$ and 
$\PP(\wedge^l V^{*})$ 
one finds

\begin{lemma}\label{gdual} $D(M)=\L(M)^{\bot}\cap G(l,m)^*$, where 
$\L(M)$ is the 
linear span of $M$.\end{lemma}
\begin{proof}  A point $P\in G(l,m)^*$ lies in $D(M)$ if and only if 
the hyperplane $H_{P}\subset \PP(\wedge^l V)$ defined by $P$ 
contains $M$, i.e. the span $\L(M)$.\end{proof}

Therefore, if $M\subset G(l,m)$ is a general subset that spans a 
linear space of 
dimension larger than that of $G(l,m)$, then the geometric dual 
$D(M)=\emptyset$. 

But for any Schubert union $S_{U}\not= G(l,m)$, the geometric dual 
$D(S_{U})$ is nonempty.  
In fact we will show that 
$D(S_{U})$ is again a Schubert union with respect to the dual flag
 $0=A^\bot_{m}\subset A^\bot_{m-1}\subset A^\bot_{m-2},..., 
A^\bot_{1}\subset V^*$.

We denote the Pl\"ucker coordinates on $\PP(\wedge^l V^*)$ with 
respect to this flag by 
$$\{ X^*_{\alpha}| \alpha=(a_{1}, a_{2},..., a_{l}), \quad 
0<a_{1}<...< a_{l}\leq m\}.$$ 
Notice that these coordinates are natural dual to the coordinates 
$X_{\alpha}$ defined on $\PP(\wedge^l V)$.  In fact
\begin{equation}
    X^{*}_{\alpha}(X_{\beta})=0\quad 
    \text{if and only if}\quad\alpha\neq {\beta} \label{pdual}
    \end{equation}

As above we may define a grid $G_{G(l,m)^*}\subset \ZZ ^l$ for these 
Pl\"ucker coordinates:
$$G_{G(l,m)^*}=\{ \alpha=(a_{1}, a_{2},..., a_{l})|1\leq a_{1}<...< 
a_{l}\leq m\}.$$
A Schubert union $S_{U^*}\subset G(l,m)^*$ with respect to this flag 
is 
associated to a $G$-grid $G_{U^*}\subset G_{G(l,m)^*}$.

\begin{defn} \label{rev}
Let the natural map
$$rev: G_{G(l,m)}\to G_{G(l,m)^*}$$
be defined as
$$(a_{1}, a_{2},..., a_{l})\mapsto (m+1-a_{l},..., m+1-a_{2}, 
m+1-a_{l}).$$
\end{defn}
\begin{rem}
{\rm Clearly the map $rev$ 
%reverses the inclusions, and 
has a natural inverse $$rev^*: G_{G(l,m)}^*\to G_{G(l,m)}.$$

It therefore sets up a duality between the two grids that we 
naturally call 
 ${\it Grid}$ ${\it duality.}$  For the relation to Poincare duality, 
see Remark \ref{Poincare}. }
%\label{grid}
\end{rem}
\begin {defn} \label{griddual}
 Let $M$ be an arbitrary subset of the 
$G_{G(l,m)}$-grid $\{(a_1,...,a_l) \in Z^l  | 1 \le a_1 < a_2 <....<
a_l \le m\}.$ Then the grid-dual of $M$ is 
$$M^{rev}= 
\{
rev(\alpha) | \alpha \in M \}\subset G_{G(l,m)^*}.$$
\end{defn}

We are now ready to define the grid dual of a Schubert union 
$U\subset$.

\begin {defn} \label{griddu}
Let $S_{U}$ be a Schubert union in $G(l,m)$ with $G$-grid $G_U$.
Then the grid dual of $S_{U}$ is the Schubert union $S_{U^*}\subset 
G(l,m)^*$ whose
$G$-grid $G_{U^{*}}$ is the grid-dual $H_{U}^{rev}$ of $H_U$.
\end{defn}
Roughly speaking this means that one finds the dual of a Schubert
union by ``turning its $H$-grid around with the map $rev$'' and use 
it as $G$-grid.
(``Turning around'' just means taking its mirror image relative to 
the level linear space $d=i_1+i_2-3=\delta/2=m-2$ if $l=2$).
The key lemma that links grid-duality to geometric duality is:

\begin{lemma}\label{ppdual} Let $M\subset G_{G(m,l)}$, and let 
$\L(M)\subset\PP(\wedge^l V)$ be the linear space defined by the 
vanishing of the Pl\"ucker coordinates 
$\{X_{\alpha}|\alpha\in G_{G(m,l)}\setminus M\}$.  Then 
$\L(M)^\bot\subset \PP(\wedge^l V^{*})$ is defined by the vanishing 
of 
the Pl\"ucker coordinates $\{X^{*}_{\alpha}|\alpha\in 
M^{rev}\}$.\end{lemma}
\begin{proof}
    This is simply the orthogonality induced by the duality of the 
two 
    basis $\{X_{\alpha}|\alpha\in G_{G(l,m)}\}$ and 
    $\{X^{*}_{\alpha}|\alpha\in G_{G(l,m)^{*}}\}$.\end{proof}

We may now state the promised result:
\begin{thm} \label{dualities}
For a Schubert union $S_U$ in $G(l,m)$ its geometric dual $D(S_U)$, 
and its grid-dual $S_{U^*}$, are equal.
\end{thm}
\begin{proof}
The homogeneous ideal $I_{S_U}$ of $S_{U}$ modulo the ideal of 
$G(l,m)$ is generated by 
$$\{X_{\beta}| \beta\in H_{U}\}.$$
By lemmas \ref{gdual} and \ref{ppdual} a linear form in this ideal 
correspond to a point $P$ on the geometric dual $D(S_{U})$ 
if and only if $P\in G(l,m)^*$ and the nonzero Pl\"ucker coordinates 
of 
$P$ lie in $H_{U}^{rev}$, 
i.e. if and only if they lie in the grid-dual  
$S_{U^*}$. \end{proof}

\begin{figure}[ht]
     \centering
     \setlength{\unitlength}{15pt}
\noindent
\begin{picture}(7,7)
  %   \thinlines
 %     \multiput(0,0)(1,0){7}{\line(0,1){7}}
 %  \multiput(0,0)(0,1){7}{\line(1,0){7}}
   \thicklines
   \put(0,0){\line(0,1){6}}
  \put(1,0){\line(0,1){6}}
 % \put(7,0){\line(0,1){7}}
%   \put(0,7){\line(1,0){6}}
   \put(0,0){{\line(1,0){1}}}
  \put(0,1){\line(1,0){2}}
   \put(0,0){{\line(1,0){1}}}
  \put(2,1){\line(0,1){5}}
 % \put(1,1){{\line(1,0){1}}}
  \put(3,2){\line(0,1){4}}
  \put(0,2){{\line(1,0){3}}}
  \put(4,3){\line(0,1){3}}
  \put(0,3){{\line(1,0){4}}}
  \put(5,4){\line(0,1){2}}
  \put(0,4){{\line(1,0){5}}}
  \put(6,5){\line(0,1){1}}
  \put(5,5){{\line(1,0){1}}}
 % \put(7,6){\line(0,1){1}}
  \put(0,6){{\line(1,0){6}}}
   \put(3,3){{\line(0,1){2}}}
  \put(0,5){{\line(1,0){6}}}
  \put(0.4,4.3){x}  
\put(3.4,4.3){x}  
\put(4.4,4.3){x}  
  \put(0.4,4.3){x}  
\put(3.4,5.3){x}  
\put(4.4,5.3){x}
\put(5.4,5.3){x}
\put(0.4,5.3){x}
\put(0.4,5.3){x}
\put(2.4,5.3){x}
\put(3.4,3.3){x}
\put(1.4,5.3){x}
  \put(1.4,4.3){x} 
 \put(2.4,4.3){x}
  \put(0.4,3.3){0}
  \put(1.4,3.3){0} 
 \put(2.4,3.3){0}
  \put(4.1,0.0){$S_{(3,5)}$}
 \put(0.4,2.3){0}
  \put(1.4,2.3){0} 
\put(2.4,2.3){0}
 \put(0.4,1.3){0}
  \put(1.4,1.3){0}
\put(0.4,0.3){0}
\end{picture}
   \caption{}
     \label{fig:u}
\end{figure}

\begin{figure}[ht]
     \centering
     \setlength{\unitlength}{15pt}
\noindent
\begin{picture}(7,7)
  %   \thinlines
 %     \multiput(0,0)(1,0){7}{\line(0,1){7}}
 %  \multiput(0,0)(0,1){7}{\line(1,0){7}}
   \thicklines
   \put(0,0){\line(0,1){6}}
  \put(1,0){\line(0,1){6}}
 % \put(7,0){\line(0,1){7}}
%   \put(0,7){\line(1,0){6}}
   \put(0,0){{\line(1,0){1}}}
  \put(0,1){\line(1,0){2}}
   \put(0,0){{\line(1,0){1}}}
  \put(2,1){\line(0,1){5}}
 % \put(1,1){{\line(1,0){1}}}
  \put(3,2){\line(0,1){4}}
  \put(0,2){{\line(1,0){3}}}
  \put(4,3){\line(0,1){3}}
  \put(0,3){{\line(1,0){4}}}
  \put(5,4){\line(0,1){2}}
  \put(0,4){{\line(1,0){5}}}
  \put(6,5){\line(0,1){1}}
  \put(5,5){{\line(1,0){1}}}
 % \put(7,6){\line(0,1){1}}
  \put(0,6){{\line(1,0){6}}}
   \put(3,3){{\line(0,1){2}}}
  \put(0,5){{\line(1,0){6}}}
%\put(0.4,4.3){x}  
\put(3.4,4.3){x}  
\put(4.4,4.3){x}  
\put(0.4,4.3){0}  
\put(3.4,5.3){x}  
\put(4.4,5.3){x}
\put(5.4,5.3){x}
%\put(0.4,5.3){x}
\put(0.4,5.3){0}
\put(2.4,5.3){x}
\put(3.4,3.3){x}
\put(1.4,5.3){0}
  \put(1.4,4.3){0} 
 \put(2.4,4.3){x}
  \put(0.4,3.3){0}
  \put(1.4,3.3){0} 
 \put(2.4,3.3){x}
  \put(4.1,0.0){$S_{(2,7)}*$$\cup S_{(3,4)}*$ }
%is the dual of $S_{(3,5)}$}
 \put(0.4,2.3){0}
  \put(1.4,2.3){0} 
\put(2.4,2.3){0}
 \put(0.4,1.3){0}
  \put(1.4,1.3){0}

\put(0.4,0.3){0}
\end{picture}
    \caption{}
     \label{fig:t}
\end{figure}

\begin{exa}
\begin{enumerate}
    \item
    Let $l=2, m=7$, and consider the Schubert cycle $S_{(3,5)}$.
    Its G-grid is 
    
$$G_{(3,5)}=\{(1,2),(1,3),(1,4),(1,5),(2,3),(2,4)(2,5),(3,4),(3,5)\},$$ 
    
and $$H_{(3,5)}=\{(1,6), 
    (1,7),(2,6),(2,7),(3,6),(3,7),(4,5),(4,6),(4,7),\ldots,(6,7)\}.$$
    Its grid-dual $S_{U^*}$ has G-grid
$$G_{(3,5)^{*}}=\{(1,2),(1,3),\ldots,(1,7),(2,3),\ldots,(2,7),(3,4)\}$$

Its geometric dual is therefore the Schubert union
$$D(S_{(3,5)})=S_{(2,7)^{*}}\cup S_{(3,4)^{*}}.$$
In Figure \ref{fig:u} we have marked  the points of the $G$-grid of 
$S_{(3,5)}$ by 0 and the points of its $H$-grid by x. In Figure 
\ref{fig:t} we have marked  the points of the $G$-grid of the dual 
union of $S_{(3,5)}$ by 0 and the points of its $H$-grid by x.
\item The geometric dual of $S_{(3,7)}$ is $S_{(3,4)^{*}}$.
\item The geometric dual of $S_{(5,6)}$ is $S_{(1,7)^{*}}$.
\end{enumerate}
\end{exa}

It would of course be nice also to have an explicit expression for 
the dual of a Schubert union.

Let $U$ be the Schubert union
\[S_{(a_{1,1}, a_{1,2},..a_{1,l})} \cup 
S_{(a_{2,1},a_{2,2},...a_{2,l})} \cup
.......\cup S_{(a_{s,1},a_{s,2},....,a_{s,l})}.\]

or phrased differently:
\[\cup_i S_{(a_{i,1},a_{i,2}....,,a_{i,l})} \]
with $ i \in I=\{1,....,s\}.$
The dual Schubert union is then given as the union over all disjoint 
partitions
$\{A_1,...,A_l\}$ of $I$ of the Schubert cycles that are described as 
follows:

\[X_1 \le m-\max\{a_{i,l}| | i \in A_1\},\]
\[X_2 \le       m-\max\{a_{i,l-1} | i \in A_2\},\]
\[              .........................\]
\[X_{l-1} \le   m-\max\{a_{i,2} | i \in A_{l-1}\},\]
\[X_l \le       m-\max\{ a_{i,1} | i \in A_l\}.\]

Each such collection of $l$ simple conditions give the Schubert union
$S_{(f_1, f_2...,f_l)}$, where

\[f_1 = \min \{f_2-1, m-\max \{a_{i,l}| i \in A_l\}\}, \]
\[f_2 = \min \{f_3-1 ,m-\max \{a_{i,l-1\}| i \in A_{l-1}}\}, \] 
\[  .........................\]
\[f_{l-1} = \min \{f_l-1, m-\max \{a_{i,2} | i \in A_2)\}\},\]
\[f_l =  m-\max \{0, a_{i,1}, i \in A_1\}.\]

This must be interpreted such that $S_{(f_1, f_2...,f_l)}=\emptyset$, 
if $f_i<i$ for some $i$.

\subsection{Further properties of dual Schubert unions} \label{techn}

Theorem \ref{dualities} obviously gives the following:
\begin{cor} \label{bidual}
The geometrically dual Schubert union of the geometrically dual 
Schubert union of the Schubert union $U$ is $U$.
\end{cor}
\begin{proof}
The analogous result for the grid dual, that is: grid-biduality,
obviously holds, and using Theorem \ref{dualities} we also have
geometric biduality.
\end{proof}
%\begin{cor} \label{fillspace}
%For a Schubert union $U$ we have $\L(U) \cap G(l,m)=U$. 
%\end{cor}
%\begin{proof}
%This follows directly from Proposition \ref{linsec}, (iv), and
%Corollary \ref{bidual}. 
%\end{proof}

It would be nice to be able to count the Schubert unions of $G(l,m)$
(relatively to a fixed flag, as always). For general $l$ and $m$ it 
is a combinatorical challenge to do this, a challenge we have not yet
met. For $l=2$, however, we will see in the next section how one 
can find the number of Schubert unions for fixed $m$.

\begin{defn} \label{dualunion}

For a Schubert union $U$ (or $S_U$) denote by $U^*$ (or $S_{U^*}$ its 
dual union.
\end{defn} 

\begin{rem} \label{part}

{\rm For the sake of completeness we will now will give a slight 
variation
of the representation of 
a Schubert union by its $G$-grid or $H$-grid.
This will not, however, bring any essentially new.
Look at a rectangular $l \times (m-l)$-box $B$ with $l$ columns and 
$m-l$ rows, formed by $l(m-l)$ squares of sidelength $1$.
We are now interested in partitions, whose Young diagrams (or Ferrers 
diagrams, see \cite{Fu}, p. 2-3 for a discussion of notation) can be 
placed
inside $B$. Obviously, these are some of the partitions of the 
integers
$N$ that range from $0$ to  $l(m-l)$.
The Young diagram of $s_1+s_2+...+s_r$, for non-increasing $s_j$,
consists of $s_i$ concecutive 
squares in row $i$, for $i=1,...,r$, where all left ends of the 
sequences of
squares start just below each other, and row $i$ is above row $j$ if 
$i<j$. 

The set of points in the $G_{G(l,m)}$-grid can in a
natural way be identified with the set  $\P_{(l,m)}$ of partitions
whose Young diagram can be placed inside $B$.
For a given partition $P$, let $c_j$ be the number of  summands $s_i$ 
in the
partition such that $s_i=j$. 
Hence, if the partition for example is $2+2+3=7$, then $c_2=2, c_3=1$ 
and $c_j=0$ for all other $j$. Obviously $P$ is characterized by the 
$c_i$.

(i) Put 
\[X_1=c_l+1, X_2=c_{l-1}+c_l+2, ...., X_{l-1}=c_2+c_3+...+c_l+l-1,
X_l=c_1+c_2+....+c_l+l.\]
Then $h(P)=(X_1,X_2,...,X_l)$ gives the natural bijection from 
 $\P_{(l,m)}$ to $G_{G(l,m)}$.

(ii) For the Schubert cycle $S_{\alpha}$ with 
$\alpha=(a_1,a_2,...,a_l)$
we have $h^{-1}(G_{\alpha})=\P_{\alpha}$, where $\P_{\alpha}$ is the 
set of
partitions $P$ with
\[c_l \le a_1-1, c_{l-1}+c_l \le a_2-2,......,c_2+c_3+...+c_l+l-1 \le
a_{l-1}-(l-1), c_1+c_2+....+c_l \le a_l-l.\]

(iii) For a Schubert union $U$ we see that $h^{-1}(G_{\alpha}$) is a
corresponding union of sets of type $\P_{\alpha}$. 

(iv) If we represents a partition $P$ by its $l$-tuple
$(c_1,...,c_l)$, then 

$\P_{(l,m)}=h^{-1}(G_{G(l,m)})=
\{(c_1,...,c_l) \in $${\bf Z}$$^l \quad | \quad c_1+...+c_l
\le m-l,$ and $c_i \ge 0,$ all $i \},$ 
in other words an alternative,
``twisted'' grid with $ m \choose l$ elements.

(v) Of course $h^{-1}(H_{U})$ is the complement of  $h^{-1}(G_{U})$
in  $h^{-1}(G_{G(l,m)})$. 
The relationship between $h^{-1}(H_{U})$ 
and $h^{-1}(G_{U^*})$ is, however, not as striking as that between 
$H_U$ and $G_{U^*}$, 
a slightly more complicated operation than just ``turning 
$h^{-1}(H_{U})$ 
around'' is necessary to obtain $h^{-1}(G_{U^*})$.

(vi) In Corollary  \ref{int-un}
we determined the number of ($F_q$-rational) points on 
$S_{U}$ as the sum of terms $q^{x_1+...+x_l-l(l+1)/2}$, where the sum
is taken over all tuples $(x_1,....,x_l)$ in the (usual)grid $G_U$.
We convert $x_i$ to $c_j$ and observe:
 \[x_1+...+x_l-l(l+1)/2=c_1+2c_2+...+lc_l.\] 
But this is the number, say $N(P)$, of which the partition $P$ in
question is a partition. Hence we see that the number of
($F_q$-rational) points on 
$S_{U}$ is the sum of terms $q^{N(P)}$, where the sum
is taken over all partitions represented by tuples  in the (usual) 
grid $G_U$.}
\end{rem}
\begin{rem}\label{Poincare}
{\rm For Schubert cycles there is of course another, well-known
  duality, Poincare duality. We know that the Poincare dual
  of $S_{(a_1,a_2,...,a_l)}$  is 
$$S_{(m+1-a_l,m+1-a_{l-1},...., m+1-a_2,m+1-a_1)}.$$ 

Using the map $rev$ defined earlier in this section this means that
the Poincare dual of $S_{\alpha}$ is $S_{rev(\alpha)}$.
Hence one can view Poincare duality as a map ($rev$) that sends 
individual points
of the $G_{G(l,m)}$-grid to their image points, 
%around the level linear space
%$d=i_1+i_2+...i_l-l(l+1)/2=\delta/2=l(m-l)/2$, 
while our duality of Schubert unions sends configurations of points 
($H$-grids) to their corresponding configurations of image points,
and in addition interchanges the roles of $G$-grids and $H$-grids.

The sum of the Krull dimensions of a Schubert cycle and its Poincare
dual cycle is $\delta=l(m-l)$, the Krull dimension of $G(l,m)$. The 
sum of
  the (affine) spanning dimension of
a Schubert cycle/union and its geometric/grid dual union is $k=$$m 
\choose
  l$, the (affine) spanning dimension of $G(l,m)$. While geometrical
  or grid duality will play an important role for us, we will not be
  concerned with Poincare duality in this paper (apart from the fact 
we have now pointed out, that the map $rev$ that we use, is 
essentially Poincare duality). 
}
\end{rem}

Let us finish this section by some remarks concerning duality of
Schubert unions over finite fields. Let $U$ be a Schubert union, and
let $g_U(q)$ be its number of $F_q$-rational points, as given
by Corollary \ref{int-un}.
Let $\delta=l(m-l)$ be the Krull dimension of $G(l,m)$.
Denote by $n(q)$ the number of $F_q$-rational points of $G(l,m)$,
and set $h_U(q)=n(q)-g(q)$.
\begin{prop} \label{dualnumbers}
Let $U^*$ denote the dual of a Schubert union $U$. Then the number of 
$F_q$-rational points of $U^*$ is $q^{\delta}h_U(q^{-1})$.
\end{prop}
\begin{proof}
This is clear since $n(q)$ is the sum of 
terms $q^{i_1+...+i_l-l(l+1)/2}$, where the sum is taken over all
points $(i_1,....,i_l)$ in the grid $G_{G(l,m)}$, and $g_U(q)$
is the corresponding sum over all points of $G_U$, and $h_U(q)$ is
the corresponding sum over all points of $H_U$. Passing from
$H_U$ to $rev(H_U)$ gives rise to the passage from  $h_U(q)$ to
$q^{\delta}h_U(q^{-1})$.
\end{proof}

\section{Schubert Unions in $G(2,m)$} 
\label{twodim}
In this section we will make a special study of Schubert unions in
$G(2,m)$ for $m \ge 3$. We will show that for each $m$ there are
$2^{m-1}$ such unions (for fixed flag) and that they in a natural way 
correspond to the set of subsets $P(M)$ of a given set $M$ with $m-1$ 
elements., we may assume $M=\{1,2,....,m-1\}$.  Moreover, taking
complements in $M$ corresponds to the duality of Schubert unions
described in general in Section \ref{dual}. We will also interpret
Schubert unions in other ways, including a more  ``physical'' one.

\subsection{Subsets of $M$ and increasing sequences of $2s$ numbers} 

We call a union of two Schubert cycles proper if it not equal to a
Schubert cycle.
Let $m=2$. We start out with the following easy, but important
 technical result:
\begin{lemma}
The union of two Schubert cycles $S_{(a,b)}$ and $S_{(c,d)}$
is proper if and only if $a < c < d < b$ or $c < a < b < d$.
\end {lemma}
\begin{proof}
The result follows from  $S_{(a,b)} \cap S_{(c,d)}=S_{(e,f)},$
where $e=\min\{a,c\}$ and $f=\min\{b,d\}.$
\end{proof}
This means that specifying a proper union of two Schubert cycles
amounts to specifying four integers $1 \le a < c < d < b \le m$,  
and then we obtain the union  $S_{(a,b)} \cup S_{(c,d)}.$
 In general, specifying a Schubert union of $s$ Schubert cycles, which
is not a union of $s-1$ Schubert cycles (a proper union of $s$
Schubert cycles), amounts to specifying
$2s$ integers $1 \le a_1 < a_2 < ...< a_s <  b_s <.... b_2 < b_1 \le
m,$ and then we obtain the union  $S_{(a_1,b_1} \cup....\cup  
S_{(a_s,b_s)}.$
This means that there is $m \choose 0$$=1$ empty set, there are $m 
\choose 2 $
unions that are just Schubert cycles, and in general there are 
$m \choose\ 2s$ unions that are proper unions of $s$ cycles, for each
$s$ up to $\frac{m}{2}$. Since 

$m \choose\ 0$$+$$m \choose\ 2$$+$$m
\choose 4$$+.... = $$m-1 \choose\ 0$$+$$m-1 \choose\
1$$+$$m-1 \choose\ 2$$+$$m-1 \choose\ 3$$+$$m-1 \choose\ 4$$+....+$
$m-1 \choose\ m-1$$= 2^{m-1}, \quad \quad \quad \quad$ 
we see that there is a potential for 
expressing the set of Schubert unions as a power set $P(M)$ as
described  above. 
Concretely, we choose to do it as follows: Think of the 
$G_{G(2,m)}$-grid as
consisting of squares with sidelengths $1$ centered at the points
of ${\bf Z}$$^2$ in the $(x,y)-$plane considered earlier. 
Then, specifying the $G$-grid of a Schubert union corresponds to
picking 
$m_1$ squares for $x=1, m_2$ squares for $x=2, ..., m_r$
squares for $x=r$, for some $r \le m-1$. Moreover $m_1 > m_2 > ... >
m_r$. This gives rise to:
\begin{defn} \label{mu}
For a Schubert union $S_U$ the (unordered) subset $M_{U}$ of 
$M=\{1,2,....,m-1\}$
is $\{m_1,m_2,...,m_r\}$ if there are exactly $m_i$ points $(x,y)$ in
$G_U$ with $x=i$, for $i=1,..,r$, and no points with $x=i$, for $i 
>r$.  
\end{defn}
If one prefers to list the numbers in increasing order, which in fact
will be strictly increasing, we have  $M_U=\{m_r,m_{r-1},...,m_1\}$.
A simple look at the plane diagram of $G_{G(2,m)}$,
represented as $m \choose 2$ squares forming a triangle, reveals that
the complement of $M_{U}$ is equal to $M_{U^*}$.

\subsection{Different descriptions of Schubert
  unions for $l=2$}

Look at the set of sequences $1 \le a_1<a_2<...a_d<b_d<...<b_1 \le m$.
Given such a sequence that represents the Schubert union 
$S_U=S_{(a_1,b_1)} \cup....\cup  S_{(a_s,b_s)},$ we give the
following:
\begin{defn} \label{sigma}
For a Schubert union $S_U$ we set $\sigma_U= 
a_1<a_2<...a_d<b_d<...<b_1 $
\end{defn}
We now will give the function $f$ that, with $\sigma_U$ as input,
gives the subset $M_U=f(\sigma_U)$ as output.
\begin{prop} \label{tranform}
If $\sigma_U$ is the sequence representing $S_U$, then
\[M_U=f((a_1<a_2<...a_s<b_s<...<b_1))=\]
\[\{b_s-a_s,
b_s-a_s+1,...,b_s-a_{s-1}-1,b_{s-1}-a_{s-1},...,b_{s-1}-a_{s-2}-1,...,b_1-a_1,...,b_1-1\}.\]

\end{prop}

We observe that the number of ``jumps'' in $M_U$ is $s-1$.
To go the opposite way we do as follows:
Given an element $M_U$ of $P(M)$, write it in increasing order as:

\[\{c_0,c_0+1,...,c_0+c_1, c_0+c_1+d_1, 
c_0+c_1+d_1+1,...,c_0+c_1+d_1+c_2,
.............,\]
\[ 
c_0+..+c_{s-1}+d_1+..+d_{s-1}+1,....,c_0+..+c_{s-1}+d_1+..+d_{s-1}+ 
c_s\}\]
where all integers $\{ d_1,...,d_{s-1}\}$ are at least two and
represent the ``jumps'' in $M_U$ (The cardinality of $M_U$ is then 
$c_1+..+c_s+s$.)
The formula $f^{-1}$ which gives the sequence $\sigma_U$ is now:
\[f^{-1}(\{c_0,..., c_0+..+c_s-1+d_1+..+d_{s-1}+c_s\})=\]
\[c_s+1<c_{s-1}+c_s+2< 
c_{s-2}+c_{s-1}+c_s+3<...<c_s+....+c_s+s-1<c_1+....c_s+s<\]
\[ c_0+c_1+....c_s+s< 
c_0+c_1+c_2+....c_s+d_1+s-1<c_0+c_1+c_2+....c_s+d_1+d_2+s-2<...<\]
\[c_0+..+c_{s-1}+c_s+d_1+..+d_{s-2}+2< 
c_0+..+c_{s-1}+c_s+d_1+..+d_{s-1}+1.\]
We leave it to the reader to check that $f(f^{-1}(M_U))=M_U$ and that 
$f^{-1}(f(\sigma_U))=\sigma_U$.

On the level of sequences $\sigma_U$ we have:
\begin{lemma} \label{dualsigma}
If $\sigma_U=\{a_1<a_2<...a_s<b_s<...<b_1\}$ then the dual sequence
$\sigma_{U^*}$ is obtained by listing the $2s+2$ integers 
$m-b_1,m-b_2,...,m-b_s, m-a_s-1, m-a_s, m-a_{s-1},...,m-a_1,m$,
with the convention that if and only if $b_1=m$, then we remove the
outer pair $m-b_1$ and $m$, and if and only if $m-a_s-1=m-a_s$, that
is: $a_s=a_{s-1}+1$, then we remove the midpair $ m-a_s-1$ and 
$m-a_s$.
\end{lemma}
\begin{proof}
Calculate $f(\sigma_U)$, take its complement $C$ in $M$, and find 
$g(C)$.
We leave the calculations to the reader.
\end{proof}
We then immediately get the following result, which also follows from 
a direct inspection of the grid diagrams involved:
\begin{cor} \label{moredualsigma}
The dual of a proper Schubert union of $s$ Schubert cycles is a proper
union of $s-1, s$ or $s+1$ Schubert cycles.  
\end{cor}

We now give a simple example, demonstrating the various
ways of representing a Scubert union and its dual.

\begin{exa} \label{manyways}
{\rm Look at $S_U = S_{(1,7)} \cup S_{(3,5)}$ in $G(2,7)$.
We see that  $G_{G(2,7)}$ roughly speaking consists of a triangular 
grid 
of integral points with corners $(1,2), (1,7),(6,7)$. One might embed
these points in squares with sidelength $1$ if one prefers.
We see that $G_U$ corresponds to $6$ squares in the first column,
$3$ squares in the second one, and $2$ squares in the third column.
Hence we have $M_U=\{2,3,6\}$. The complement of $M_U$ in $M$ is
$\{1,4,5\}.$ We observe that both these sets have one ``jump'' each
(from $3$ to $6$ and $1$ to $4$, respectively), and hence each of them
correspond to a proper union of exactly two Schubert cycles.
Furthermore we see that $H_U$ consists of $5$ squares in 
the upper row, $4$ squares in the row below, and $1$ square in the
row that is third from the top. ``Turning this around'' (using the
operation $(a,b) \khpil\ (8,8)-(b,a)$ from Definition \ref{griddual}) 
we get
$5$ squares in column $1$ and $4$ squares in column $2$, and $1$
square in column $3$, in other words we get $G_{U^*}$.
We observe that  $G_{U^*}$ is the union of $S_{(2,6)}$ and 
$S_{(3,4)}$.
Hence $S_{U^*}=S_{(2,6)} \cup S_{(3,4)}$.

Moreover we see that $\sigma_U= \{1 < 3< 5 < 7\}$, while 
$\sigma_{U^*}= \{2 < 3 <4 < 6\}.$ We encourage the reader to 
check the various duality formulas and transformations $f$ and $g$
described above in this example. 

%Using the terminology of Remark \ref{part}, and the ``alternative
%grid'' introduced there, we see that $S_{(1,7)}$
%corresponds to $\{(c_1,0) | 0 \le c_1 \le 5\}$, while 
%$S_{(3,5)}$ corresponds to $\{(c_1,0) | 0 \le c_1 \le 3\} \cup 
% \{(c_1,1) | 0 \le c_1 \le 2\} \cup  \{(c_1,2) | 0 \le c_1 \le 2\},$
%and the union $S_U$ thus corresponds to the union of these $2$ sets,
%that is: a set of $6$ grid points for $c_2=0$, $3$ grid points for
%$c_2=1$, and $2$ grid points for $c_2=2$.

Using the terminology of Remark \ref{part}, and the ``alternative
grid'' introduced there, we see that $S_{(1,7)}$
corresponds to $\{(c_1,0) | 0 \le c_1 \le 5\}$, while 
$S_{(3,5)}$ corresponds to $\{(c_1,0) | 0 \le c_1 \le 3\} \cup 
 \{(c_1,1) | 0 \le c_1 \le 2\} \cup  \{(c_1,2) | 0 \le c_1 \le 2\},$
and the union $S_U$ thus corresponds to the union of these $2$ sets,
that is: a set of $6$ grid points for $c_2=0$, $3$ grid points for
$c_2=1$, and $2$ grid points for $c_2=2$.
Likewise the dual union $S_{U^*}$ corresponds to $5$ gridpoints for
$c_2=0, 4$ grid points for $c_2=1,$ and one grid point for $c_2=2$. 
Likewise the dual union $S_{U^*}$ corresponds to $5$ gridpoints for
$c_2=0, 4$ grid points for $c_2=1,$ and one grid point for $c_2=2$.
}
 
\end{exa}

We end this section with a more ``physical'' remark. 
\begin{rem} \label{friction}
{\rm We represent
the $G_{G(2,m)}$ by $m \choose\ 2$ squares in a triangle, as
described. Rotate this configuration of squares an angle
$\frac{\pi}{4}$ counterclockwise, and upscale by a factor $\sqrt{2}$ 
in each direction, so that the point
$(1,2)$ stays fixed, the corner point $(m-1,m)$ is moved to
$(1,2m-4)$, and the third corner point $(1,m)$ is moved to $(3-m,m)$.
Assume now that you start with any collection of $K \le$$m \choose\ 2$
squares that you put inside the ``frame'' $G_{G(2,m)}$.
After we have rotated the frame, we let vertical gravity work, and we 
assume
that the $K$ squares can move in a frictionless way inside the frame
and relatively to each other.
Then the configuration of $K$ squares stays in equilibrium if and only
of they form  a $G_U$-grid for a Schubert union
$S_U$. The function $x+y-3$ before, and simply
$(y-2)$ after, rotating, can be thought of as a potential. The 
highest peak(s)
of $G_U$ has(have) height equal to the Krull dimension of $S_U$, 
which 
by definition is the highest Krull dimension among the Schubert
cycles, of which $S_U$ is a union.
If we let inverse gravity work, then
a configuration stays in equilibrium if and only if it forms a 
$H_U$-grid.
 
An obvious analogous formulation can be given for $l=3$, and using
slightly less concrete formulations, also for $l \ge 4$.}
\end{rem}

\section{Applications to Codes}
\label{codes}

Let $C$ be the Grassmann code $C(l,m)$ over a finite field $F_q$ as 
described for example in \cite{N}, \cite{GL}, and \cite{GT}.
Here $l$ and $m$ are natural numbers with $l < m$.
The code  $C(l,m)$ is defined as follows. 
One starts with $G(l,m)$, which is embedded in the Pl\"ucker space
$\PP^{k-1}$ with $k=$$m \choose\ l$.
It is well known that $G(l,m)$ contains $n$ points, where 

\begin{equation} \label{number}
n=\frac{(q^m-1)(q^{m-1}-1) \dots 
(q^{m-l+1}-1)}{(q^l-1)(q^{l-1}-1)\dots (q-1)}.
\end{equation}
For the special case $l=2$, which we will often study, this is 
clearly:
\begin{equation} \label{number2}
n=\frac{(q^m-1)(q^{m-1}-1)}{(q^2-1)(q-1)}.
\end{equation}

Of course both these formulas are special cases of the last point of 
Corollary \ref{int-un}. Pick a representative of each of the $n$ 
points 
as a column vector in $(F_q)^k$, and form a $k \times n$-matrix $M$ 
with 
these $n$  vectors
as columns (in any preferred order). The code  $C(l,m)$ is then the
code with $M$ as generator matrix. Hence $C$ is a linear $[n,k]$-code
(only defined up to code equivalence, since we have not specified
which representative in $(F_q)^k$ we choose for each point, but this
ambiguity will not affect the questions we will study, concerning code
parameters and higher weights).

We will now recall and establish some facts about the weight 
hierarchy  
$d_1 < d_2 <....< d_k$
of the codes $C(l,m)$, using what we know from the previous sections 
about Schubert unions. 
%More importantly, we believe, we will raise some questions. 

It has been shown in \cite{N} that the higher weights $d_r$ satisfy
\begin{prop} \label{Nogin}
$d_r=q^{\delta}+q^{\delta -1}+ \dots +q^{\delta -r+1}$, for 
$r= 1, \dots s$, where $s= max(l,m-l)+1$, and $\delta = \dim{G(l,m)}
=l(m-l).$
\end{prop} 
Now $s$ is in almost all cases much smaller than $k$, so there still
remains a lot to be shown. The proof in \cite{N} involves $3$ 
elements, the
first one is a special proof, using mulitilinear algebra, for 
$d_1=q^{\delta}$. The second ingredient is the socalled Griesmer
bound,
valid for all linear codes, in our case it gives
(for all $r$ in the range $2,...,k$, and in particular $2,....,s$):
 $$d_r \ge \sum_0^{r-1}\frac{d_1}{q^i}. $$
This gives:
\[d_r \ge q^{\delta}+q^{\delta -1}+ \dots +q^{\delta -r+1},\]
for $r=1,...,s$.
The third ingredient is the usage of the following well known fact:
Let $S$ be the set of column vectors in a generator matrix for a
linear code $C$. Then its higher weights satisfy:
\begin{equation} \label{r-weight}
  d_r = n-H_r,  
\end{equation}
for all $r=1,...,\dim C$,
where $H_r$ is the maximum number of points from $S$ contained in
a codimension $r$ subspace of $(F_q)^k$. In our case the columns are
the points of $G(l,m)$, so once and for all we define:
\begin{defn} \label{Hr}
$H_r$ (or if necessary to specify $l$ and $m$: $H_r^{l,m}$) is the
maximum number of points from $G(l,m)$ contained in a codimension $r$ 
subspace of $\PP^{k-1}$.
\end{defn}
In \cite{N} one  exhibits concrete codimension $r$ 
subspaces of $\PP^{k-1}$ containing $n -(q^{\delta}+q^{\delta -1}+
\dots +q^{\delta -r+1})$ points, for $r=1,...,s$. Equation 
\ref{r-weight} then gives:
\[d_r \le q^{\delta}+q^{\delta -1}+ \dots +q^{\delta -r+1},\]
which in conjunction with the Griesmer bound gives the result.
The third ingredient, or step, has been given in an alternative way in
\cite{GL},
using socalled close families (of what we would call grid points).
In \cite{GL} one looks at the zero set within $G(l,m)$
of $X_{\beta_1},....,X_{\beta_r}$, where the $\beta_{i}$ are grid
points with $l-1$ common coordinates (this is what it means that the 
family is close). One then shows that this zero
set contains $n -(q^{\delta}+q^{\delta -1}+
\dots +q^{\delta -r+1})$ points, in analogy with \cite{N}, and that 
you
can find close families up for $r=1,...,s$.  

We will now use Schubert unions to give a third (admittedly similar)
variant of this third step.
We will, however, first show a dual result:
\begin{thm} \label{highweights}
For the $q$-ary code $C(l,m)$ defined in the introduction
we have $d_k=n$, and
\[d_{k-a}=n-(1+q+ \dots q^{a-1}),\]
for $a=1, \dots ,s$.
 \end{thm}
\begin{proof}
Without loss of generality we assume $s=m-l$. 
We use the notation and flag $A_0 < A_1 <..... < A_m$ defined in the
start of Section \ref{SUdef}. Look at the set of 
$l$-dimensional subspaces of $V$ that contain $A_{l-1}$. This behaves
under the Pl\"ucker embedding as the set of lines through the origin
in $V'=(F_q)^{m-l+1}$, in other words it is a projective space of
dimension $s-1=m-l$. This space then of course contains projective
subspaces of dimension $i$, for $i=0,....,m-l$. In particular we look
at the set $S$ of those $l$-spaces $V$ that contain  $A_{l-1}$ and are
contained in $A_{l+i}$, for $i=0,...,m-l$. But this is by definition
the Schubert cycle $S_{\alpha_i}$, where $\alpha_i=(1,2,...,l-1,l+i)$.
In particular it is a Schubert union $U_i$ with $G$-grid 
\[G_i=\{(1,2,..,l-1,l), (1,2,..,l-1,l+1, (1,2,..,l-1,l+2),....,
(1,2,...,l-1,l+i)\}.\]
It contains $1+q+q^2+.....+q^i$ points in virtue of being a 
projective space
or alternatively, by using Proposition \ref{Sp} (or Corollary 
\ref{int-un}).
This is a priori the maximum number of points (from $G(l,m)$ or even
from all of the Pl\"ucker space) a $(k-1-i)$-codimensional subspace 
of the Pl\"ucker
space may contain. This means that these Schubert unions compute
$H_r$, and therefore $d_r$, for $r=k-1-m+l,....,k-1$, in other words
for $m-l=s-1$ consecutive values of $r$.This gives the result.  In 
addition we have,
trivially, $d_k=n$, computed by the particular Schubert union
$\emptyset$. 
\end{proof}
\begin{rem} \label{dualproof}
{\rm Theorem \ref{highweights} uses only the well-known fact of 
the socalled index of the Grassmann variety, the maximum dimension of 
linear subspace of $G(l,m)$. We chose to present a detailed proof 
here 
in order to
demonstrate how in some cases $H_r$ is computed by Schubert unions, 
and 
to give an easy variant of Step $3$ of the proofs in \cite{N} and
\cite{GL} of Proposition \ref{Nogin}, as follows: 

For each $i=1,...,m-l$, we study 
%use Proposition \ref{dualnumbers} to
%find the number of points on 
the Schubert unions $U_i^*$ that are dual to the
ones with $G$-grid $G_i$ as just described.
that is the ones with $H$-grid equal to $G_i^{rev}$.
Set $\delta=l(m-l)$.
 Proposition \ref{dualnumbers} now gives that the number of points on 
$U_i$ is $$q^{l(m-l)}(n(q^{-1})-(1+q^{-1}+.....+q^{-i})),$$ for 
$i=1,...,m-l$, and this is equal to:
$$n(q)-(q^{\delta} + q^{\delta-1}+.....+q^{\delta-i}).$$
Since the affine spanning dimension of $U_i$ is $i+1$ for each $i$,
the spanning dimension of $U_i^*$ is $k-(i+1)$, and hence, by 
Corollary \ref{int-un}, $U_i^*$ is in fact equal to a linear section 
of
$G(l,m)$ of codimension $r=i+1$. As a consequence we
obtain:
$H_r \ge n(q)-(q^{\delta} + q^{\delta-1}+.....+q^{\delta-(r-1)})$,
for $r=2,...,m-l+1=s$, and therefore:
$$d_r \le q^{\delta} + q^{\delta-1}+.....+q^{\delta-(r-1)}.$$
This gives an alternative proof of Step $3$.
Even if all the $U_i$ are Schubert cycles, many of the $U_i^*$
are proper unions of more than one cycle, in codimension $2$, for
example, one needs more than one cycle of Krull dimension
$\delta-2$.  
One sees that the points of the $H$-grids of the $U_i^*$ satisfy the
requirements of being a close family, in the sense of \cite{GL}, 
since 
the points of $G_i$ do, and the $H$-grids in question are obtained
by turning around the $G_i$ as described (the only coordinate that
varies for the points in these $H$-grids will be the first one, since
the only one varying in each $G_i$ is the last one). Hence one could
use the results in \cite{GL} instead of  Proposition \ref{dualnumbers}
to find the number of points on the $U_i^*$ also.}
 \end{rem}

\begin{example} \label{CIproj}
{\rm For $C(2,4)$ we have $k=6$ and $s=3$. Nogin's result gives
$d_1=q^4, d_2=q^4+q^3, d_3=q^4+q^3+q^2$. Theorem  \ref{highweights}
gives $d_6=n, d_5=n-1, d_4=n-1-q, d_3=n-1-q-q^2$. Comparing the two
formulas for $d_3$, we obtain $n=q^4+q^3+2q^2+q+1$, which  
is equal to the well-known formula 
$\frac{(q^4-1)(q^3-1)}{(q^2-1)(q-1)}$.
We remark that the Griesmer bound gives: 
$d_4 \ge q^4+q^3+q^2+q=d_3+q$, while the true value is
$q^4+q^3+2q^2=d_3+q^2$.
This code has been studied in much greater detail in \cite{MV}
where one has not only calculated the higher weights, but also the
higher spectra (how many subspaces of dimension $r$ of the given code
$C(2,4)$ have support weight $s$ for each conceivable $r$ and $s$) .

For $C(2,5)$ we have $k=10$ and $s=4$. Nogin's result gives
$d_1=q^6, d_2=q^6+q^5, d_3=q^6+q^5+q^4, d_4=q^6+q^5+q^4+q^3$. 
Theorem  \ref{highweights} gives 
$d_{10}=n, d_9=n-1, d_8=n-1-q, d_7=n-1-q-q^2, d_6=n-1-q-q^2-q^3 $. 
Here $n=q^6+q^5+2q^4+2q^3+2q^2+q+1$.
The only remaining case is $d_5$. The Griesmer bound gives:
$d_5 \ge d_4 + q^2$, while we know that $d_6= d_4 + q^4 + q^2$.
We have:}

\end{example}

\begin{prop} \label{d5}
For the code $C(2,5)$ we have $d_5=n-(q^3+2q^2+q+1)=d_4+q^4=d_6-q^2.$ 
\end{prop} 

\begin{proof}
The Grassmannian $G=G(2,5)$ is embedded by the 
Pl\"ucker embedding in $\PP^{9}$. Let $H_{5}$ be the maximal number 
of points in the intersection of $G$ with a  $4$-space (codimension 
$5$ ) in the Pl\"ucker $9$-space $P$. 
We use the formula $d_5=n-H_5$ to prove 
%   \[d_4+q^2=d_6-q^4 \le  d_5\] 
 \[d_5  \le  d_4+q^4=d_6-q^2\]
since  $S_{(1,5)} \cup S_{2,3}$
contains $q^3+2q^2+q+1$ points, and spans a codimension $5$ space by 
Proposition \ref{int-un}. 
%This information is summarized in 
 %  \[q^3+2q^2+q+1 \le  H_5  \le  q^4+q^3+q^2+q+1,\]
%where $H_{5}$ is the maximal number of points in the intersection of
%$G$ with a codimension $5$-space, i.e. a $4$-space in Pl\"ucker $9$-space $P$. 
To complete the proof we now prove the supplementary inequality 
$d_5  \ge d_4+q^4 = d_6-q^2$, or  equivalently:
 $H_{5} \le q^3+2q^2+q+1$. We know that
$G$ is cut out by the Pl\"ucker quadrics $Q_1,.., Q_5$ in $P$ 
since there are  in general $m  \choose\ 4$ such relations for 
$G(2,m)$, for 
$m \ge 3$.
A $4$-space $W$ in $P$ cannot be contained in $G$,  
since the maximum number $s$ such that $G$ contains an
$s$-space, is $3$.  In fact, any line in $G$ is formed by all 
projective lines through a point in 
a plane in $\PP^{4}$, so a linear subspace of $G$ is formed by either 
all  lines in a plane, 
or by all lines through a point in a linear subspace of $\PP^{4}$.  
Therefore the maximal dimension is obtained by the family of 
all lines through a point in $\PP^{4}$, which of course is 
$3$-dimensional.

Hence at least one of the Pl\"ucker quadrics does not contain $W$,
and so $W_G=W \cap G$ is contained in a quadric in $W=\PP^{4}$. 
The maximal rank of the restriction of the quadrics $Q_{i}$ to $W$ is 
between $1$ and $5$. 
In case the maximal rank is $1$ or $2$, each quadric decomposes in 
linear factors when restricted to $W$,
 so the intersection $W_{G}$ is a union of linear subspaces.  As 
above the maximal dimension of a linear component is $3$, and in that 
case it corresponds to the lines through a fixed point $P$.  
 If this is one of the components then the residual is also linear, 
and from the description above we conclude in the following way that 
it has dimension at most $2$, and the two components intersect in a 
line: If the residual component had dimension $3$, then that would 
also correspond to the lines through (another)
fixed point $Q$. On one hand the intersection between the 
$3$-dimensional components is a plane, since we are in $\PP^{4}$. On 
the other hand it only consists of one point, the one corresponding 
to the line between the two fixed points $P$ and $Q$. This is a 
contradiction, and hence the residual component is at most a plane, 
and it intersects the three-dimensional component in codimension one 
in the residual component.
 In this case the cardinality of $W_G$ is at most $q^3+2q^2+q+1$.  If 
the dimension 
 of the linear components of $W_{G}$ are smaller than $3$ and $2$ 
then $W_{G}$ is always contained in the
  union of a $3$-space and a $2$-space so the cardinality 
 is always smaller than $q^3+2q^2+q+1$.

 In case the maximal rank of the restriction of the quadrics $Q_{i}$ 
to $W$ is at least $3$, then we get the desired upper bound from the 
cardinality
of points on an irreducible quadric. By projecting to $\PP^{4}$ 
from a
smooth point one gets  maximal cardinality $q^3+q^2+q+1$ when the
quadric has rank $3$ or $5$, while the maximal cardinality is 
$q^3+2q^2+q+1$ when
the quadric has rank $4$ (cf. \cite{HT},  p 4-5 for details). 

%the cardinality 
 %of a quadric $Q$ in $W$.    But this is easily checked by the rationality of the quadric: 
 %Let $p\in Q$ be a smooth point ( if there are none, then the cardinality is at most $q+1$, the maximal cardinality of the singular points).
 %The projection of $Q$ from $p$ is one to one onto $\PP^3$ outside the lines in $Q$ through $p$.
 %The image of the lines in $Q$ through $p$ are points on a conic section $C_{p}$ in a plane. 
 %The conic has rank $2$ less than the rank of $Q$.  The cardinality of $Q$ is now counted as the cardinality of $\PP^3$ minus the cardinaliy of the plane of $C_{p}$ 
 %pluss the cardinality of the preimage of $C_{p}$.
 %In case the conic $C_{p}$ has rank $1$ it is a double line and there is a unique plane in $Q$ through $p$.  
 %¤The cardinality of $Q$ is therefore in this case $q^3+q^2+q+1-(1+q+q^2)+(1+q+q^2)=q^3+q^2+q+1$.
 %If the conic $C_{p}$ has rank $2$, there are precisely two planes in $Q$ intersecting in a line through $p$, and the cardinality of $Q$ is 
 % $q^3+q^2+q+1-(1+q+q^2)+2(1+q+q^2)-(1+q)=q^3+2q^2+q+1$.
 %  If the conic $C_{p}$ has rank $3$ there is a line in $Q$ passing though any point on $C_{p}$ and $p$, and the cardinality of $Q$ is 
 % $q^3+q^2+q+1-(1+q)+q(1+q)+1=q^3+2q^2+q+1$.  Thus the maximal cardinality of $W_{G}$ is $q^3+2q^2+q+1$, which is what we wanted to prove.
 
\end{proof}
\begin{rem}
{\rm We have observed that  $d_5$ is computed by fourspaces spanned by
Schubert unions $S_{(1,5)} \cup S_{(2,3)}$
correponding to pairs $(a,H)$, with $a$ in $H$ in $\PP^{4}$. In
addition, using Proposition \ref{int-un}, 
we see that $d_5$ is computed by 
Schubert unions (in fact cycles) of type
$S_{(2,4)}$, which happen to be duals of the $S_{(1,5)} \cup 
S_{(2,3)}$
To be concrete: For  each $3$-space $F$ in $\PP^{4}$, 
and  each line L in $F$, we look at the set of lines contained in $F$,
intersecting $L$.
%The number of lines in $F$
%is the cardinality of $G(2,4)$, namely $q^4+q^3+2q^2+q+1$.
All in all, for each pair $(a,H)$ with $a \in H$, and each pair 
$(L,F)$,
with $L \in F$, we get a $4$-space $W$ that computes $d_5$.
There are $(q^4+q^3+2q^2+q+1)(q^3+2q^2+q+1)$ of both kinds, and we
conjecture that no other $4$-spaces compute $d_5$, and that the number
of $5$-subspaces of $C(2,5)$ with minimal support weight $d_5$ 
therefore is  $2(q^4+q^3+2q^2+q+1)(q^3+2q^2+q+1)$. }

\end{rem}
\begin{defn} \label{Sch}
For given $l,m$, set $\Delta_r=d_r-d_{r-1}$ for $r=1,....,k$. 
($\Delta_0=0.$)
\end{defn}
Example \ref{CIproj} and Proposition \ref{d5} enable us to make a 
complete table for all the 
true values of $\Delta_r$ for $(l,m)= (2,5)$.
Although only $\Delta_5, \Delta_6$ for $m=5$ gives something which 
could
not be concluded from 
Proposition \ref{Nogin} and Theorem \ref{highweights}, we also, to
illustrate, give corresponding tables for $(l,m)=(2,3)$ and $(2,4)$.:

\[C(2,3):  \left[
\begin{array}{cccc} 
r:         & 1   & 2   & 3     \\
\Delta_r:  & q^2 & q & 1   

\end{array}
\right]\]

\[C(2,4):  \left[
\begin{array}{ccccccc} 
r:         & 1   & 2   & 3   & 4   & 5   & 6   \\
\Delta_r:  & q^4 & q^3 & q^2 & q^2 & q   &  1   
 
\end{array}
\right]\]

\[C(2,5):  \left[
\begin{array}{ccccccccccc} 
r:         & 1   & 2   & 3   & 4   & 5   & 6   & 7   & 8   & 9   & 10 
\\
\Delta_r:  & q^6 & q^5 & q^4 & q^3 & q^4 & q^2 & q^3 & q^2 & q   & 
1   

\end{array}
\right]\]
All values of $d_r$ in these cases are computed by Schubert unions.

This motivates the following definitions:
\begin{defn} \label{Sch2}
For given $l,m$, let $J_r$ be the maximum number of points in a
Schubert union spanning
a linear space of codimension at least $r$ in the Pl\"ucker space, 
and  
set $D_r=n-J_r$, and $E_r=D_r-D_{r-1}$, for $r=1,....,k$. ($D_0=0.)$ 
\end{defn}

We end this section  with  the following result:
\begin{prop} \label{SUbound}
For all $l,m,$ and $r$ we have
\[d_r \le D_r.\]
\end{prop}
The fact that the result is obvious does not prevent that it is 
useful,
since we have Schubert unions for any spanning
dimension, and since we can calculate the numbers of points contained
in them, using Corollary \ref{int-un}. In Section \ref{algo} we will 
give general methods to calculate the upper bound $D_r.$
Nevertheless, it is an open question whether the upper bound $D_r$ is 
equal to 
the true value $d_r$ in the cases not determined by Proposition 
\ref{Nogin}, Theorem \ref{highweights} and Preoposition \ref{d5}.

\section{Schubert unions with a maximal
 number of points}
\label{algo}

In this section we will make significant steps toward finding the 
$J_r$ and 
$D_r$ for all $m$ when $l=2$.

\begin{defn} \label{orders} 
 Fix a dimension $0\leq K\leq {m\choose l}$, and
consider the set of Schubert unions $\{U\}_{K}$ in $G(l,m)$ with 
spanning dimension $K$.

Then we order the elements $U$ in $\{U\}_{K}$ according to the 
lexicographic 
order on the polynomials $g_U$. In other words $U > V$ if $\deg g_U > 
\deg g_V$ or $\deg g_U = \deg g_V$, and the coefficient of $g_U$ is 
larger than that of $g_V$ in the largest degree where the 
coefficients differ. 
We call this the order with respect to $g_U$

For $l=2$ the elements in $\{U\}_{K}$ can also be ordered in another 
way. 
Recall Definition \ref{mu} of  $M_U=\{m_r,...m_2,m_1\}$ with $m_r 
<...<m_2<m_1$ for each $U$. 
Then $U > V$ if the largest element $m_1$ of $U$ is larger than that 
of $V$, or if these elements are equal, if the second largest, 
$m_{2}$, 
is larger than that of $V$, and so on. We call this the order with 
respect 
to $M_U$. 
\end{defn}

\begin{rem} \label{both}
{\rm It is clear that for given $m$ and $K$, we obtain the $G$-grid 
of the maximal 
element of $\{U\}_{K}$ with respect to $M_U$ by ``filling up as many 
columns of the 
$G_{Gl(2,m)}$-grid as we can from the left''. Likewise we obtain the 
$G$-grid 
of the minimal element of $\{U\}_{K}$ with respect to $M_U$ by 
``filling up as many rows of the $G_{Gl(2,m)}$-grid as we can from 
the bottom''.  

 Each  non-empty Schubert union, maximal with respect to maximal 
$M_U$, is a
 unions of
  two Schubert cycles as follows:

$S_{(x,m)} \cup S_{(x+1,y)}$, with $1 \le x \le m-1$ and $1 \le y \le 
m$.

The non-empty $U$ minimal with respect to $M_U$ are unions of type

$S_{(x,x+1)} \cup S_{(a,x+2)}$, with $1 \le x \le m-1$ and $1 \le a 
\le x+1$.}
\end{rem}
We have:
\begin{prop} \label{leftright}  Assume $l=2$. 
Fix a dimension $0\leq K\leq {m\choose 2}$, and
consider the set of Schubert unions $\{U\}_{K}$ with spanning 
dimension $K$.
Let $U_{1}$ and $U_{2}$ be the maximal and minimal elements in 
$\{U\}_{K}$ 
with respect to $M_U$.
Then $U_{1}$ or $U_{2}$ is maximal in  $\{U\}_{K}$  with respect 
to $g_U$. 
Furthermore, the one(s) that is(are) maximal with respect to $g_U$, 
also has(have) the maximum number of points over $F_q$ for all large 
enough $q$.  
\end{prop}

\begin{proof}
 Given the spanning dimension $K$, let $d=d(K)$ be the
maximal Krull dimension for the Schubert unions $\{U\}_{K}$. 
This Krull dimension is the crucial ingredient in 
our argument, since the Krull dimension is the degree of the 
polynomial $g_{U}$. We will find the maximal polynomial $g_{U}$ 
in the lexicographic order.
The fact that the union(s) that is(are) maximal with respect to 
$g_U$, also has(have) the maximum number of points over $F_q$ for all 
large enough $q$, is obvious.
Our argument is visualized by the $G_{G(2,m)}$-grid, arranged
as a set of squares in a triangle defined by
$$G_{G(2,m)}=\{(x,y)| 1\leq x<y\leq m\}.$$
Each point $(a,b)\in G_{G(2,m)}$ defines a Schubert cycle $S_{(a,b)}$ 
with Krull dimension $d(a,b)=a+b-3$.  Therefore the Schubert cycles 
with a 
fixed Krull dimension lie on the diagonal 
$$D_{d}=\{(x,y) |\quad 1\leq x<y\leq m, \quad x+y-3=d\}.$$
Let as above $$G_{a,b}=\{(x,y)\in G_{G(2,m)}| x\leq a, y\leq b\},$$
and $$G_{U}=\cup_{S_{(a,b)}\subset U}G_{a,b}.$$

By definition of $d=d(K)$, there is a
Schubert union $U$ of spanning dimension $K$ with a $G_U$ that 
contains 
a point $(a,b)$ on the diagonal $D_{d}$, i.e. $a+b-3=d$, but there is 
no such union 
with $G$-grid that contains a point on the diagonal  $D_{d+1}$. 

The {\it cardinality} $c(x,y)$ of a $G$-grid $G_{x,y}$ defines the 
function
$$c:G_{G(2,m)}\to Z, \quad (x,y)\mapsto xy- \frac{x(x+1)}{2}.$$

The restriction of this function to the diagonal $D_{d}$ is defined by
$$c(x,d-x+3)=x(d-x+3)-\frac{x(x+1)}{2}, \quad {\rm for}\quad {\rm  
max} \{d+2-m,0\} < x < \frac{d+3}{2}$$
which is clearly quadratic and concave.  Therefore it attains its 
minimum $C(d)$, when $x$ is minimal or maximal, i.e. at one of the 
end points of 
the diagonal $D_{d}$.

  Clearly $$C(d(K)) \le K \le C(d(K)+1)-1.$$ 

 We say that a point $(a,b)\in D_{d(K)}$ is {\bf admissible}, if 
$G_{a,b}\subset G_{U}$ for 
  some Schubert union $U$ of spanning dimension $K$.  Equivalently, 
  $(a,b)\in D_{d}$ is admissible if $$c(a,b)< C(d+1),$$
  i.e. has less cardinality than any point in the next diagonal.
 
 Next, we characterize the admissible points by which diagonal 
$D_{d}$ 
 they belong to.

\begin{lemma} \label{admiss}
    Consider the diagonal 
     $$D_{d}=\{(x,y)|x+y-3=d\}=\{(x,d-x+3)\quad |\quad {\rm 
max}\{d+2-m,0\}<x<\frac{d+3}{2}\}$$
    (i) Let $d\leq m-3$, then the only admissible point on $D_{d}$ is 
    $(1,d+2)$, except when $d=2$, since the point $(2,3)$ is also 
    admissible.  
   
(ii) Let $d>m-3$, then $(x,d-x+3)$ is an admissible point on the 
diagonal $D_{d}$
only if $d+3-m\leq x\leq d+4-m$ or $\frac{d}{2}\leq 
x\leq\frac{d+2}{2}$, i.e. only 
if it is among the two points with the smallest value of $x$, or the 
two 
points with largest value $x$, with one exception, namely when $m=11$ 
and $d=10$, then the point $(4,9)$ is also admissible.

(iii) If  $m>10$, then  the point $(x,m)$ is admissible, only if $x 
\ge m-3$ or $x \le
\frac{m}{5} + 2$ if $x+m$ is odd, and only if $x \ge m-3$  or  $x \le
\frac{m}{5} + 1$ if $x+m$ is even.

(iv) If  $m>10$, then  the point $(x,m-1)$ is admissible if $x \ge 
m-4$ or
$x \le \frac{m}{5} + 1$ if $x+m$ is odd, and only if $x \ge m-4$ or 
$x \le
\frac{m}{5} + 2$ if $x+m$ is even.
\end{lemma}
\begin{rem} \label{interior}
{\rm 
(a) The lemma implies, in very rough terms, that if you cross out the 
leftmost column, the two
uppermost rows, and the two right-lowest points on each diagonal
(it's really only necessary to cross out the right half of them, in
addition to $(2,3)$) then the ``interior'' grid that
remains contains no admissable points (except $(4,9)$ for $m=11$).

(b) If $d$ is odd, then part (ii) of the lemma implies that the point 
of the diagonal $D_{d}$ with the next to largest $x$, is 
non-admissible, so that the admissible ones must be among the two to 
the left, and the one to the right.}
\end{rem}
\begin{proof}
(i) This follows from an easy application of the concavity of the 
cardinality 
function $c(x,y)$ along the diagonals.  For the smaller values of $x$ 
we get 
$$c(2,d+1)=2d-1>c(1,d+3)=d+1\quad d\geq 2$$
except when $d=2$, while for the larger values of $x$ we have 
$c(\frac{d+2}{2},\frac{d+4}{2})=\frac{d^2+6d+8}{8}>d+1$ when $d$ is 
even, and 
$c(\frac{d+1}{2},\frac{d+5}{2})=\frac{d^2+8d+7}{8}>d+1$ when $d$ is 
odd.

(ii) Again, by the concavity of the function $c(x,y)$ restricted to 
the diagonal $D_{d}$, 
it suffices to compare a 
few points on each diagonal with the end points on the next.  
Furthermore the diagonal $D_{d}$ has length at least $5$ only if 
$\frac{m+1}{2}>5$ or $\frac{d+2}{2}-(d-m+3)\geq 5$, i.e. when $m>9$ 
or $d<2m-14$. So from here on we restrict to these values of $d$ and 
$m$.
Technically there is a difference between $d$ odd and even.
More precisely, if $d$ is even, we compare the value of $c$ at the 
endpoints $(d-m+4,m)$ and 
$(\frac{d+2}{2},\frac{d+6}{2})$  
on $D_{d+1}$, with its value on $(d-m+5,m-2)$ and 
$(\frac{d-2}{2},\frac{d+8}{2})$, namely the points number
 three from the endpoints on $D_{d}$.  
 Thus $(d-m+5,m-2)$ is admissible only if
 $$c(d-m+5,m-2)<c(d-m+4,m)\quad {\rm and}\quad 
c(d-m+5,m-2)<c(\frac{d+2}{2},\frac{d+6}{2}).$$
 The first inequality reduces to
 $4m-15<3d$, while the second one yields 
 $$5d^2+(70-16m)d+12m^2-100m+216>0,$$
 hence
 $$(2m-10-d)(5d-6m+20)<16.$$
 In this latter inequality, either the factors have different sign or 
 both factors have small positive value.  On the one hand 
 $d<2m-14$, means $2m-10-d>0$, so only the second factor can be 
 negative.  In this case,
 $4m-15<3d$ and $5d-6m+20<0$, which means $m\leq 7$, contrary to the 
 above.   Thus both factors have small value, i.e. $2m-26\leq d\leq 
 2m-14$ and $\frac{6}{5}m-4\leq d \leq \frac{6}{5}m-\frac{4}{5}$. 
 These inequalities are both satisfied only if 
 $2m-26<\frac{6}{5}m-\frac{4}{5}$, i.e. $m<32$.
 It is easily checked that the assertions hold for these small values 
 of $m$ except $m=11$ 
and $d=10$, when the point $(4,9)$ is also admissible. 

The point $(\frac{d-2}{2},\frac{d+8}{2})$ is admissible only if
$$c(\frac{d-2}{2},\frac{d+8}{2})<c(d-m+4,m)\quad {\rm and} \quad 
c(\frac{d-2}{2},\frac{d+8}{2})<c(\frac{d+2}{2},\frac{d+6}{2}).$$
The latter inequality yields $d<12$, and hence $m<d+3=15$, in which 
case the first inequality is satisfied only if  $m=11$ 
and $d=10$, when the point $(4,9)$ is also admissible.

If $d$ is odd, we compare the value of $c$ at the 
 endpoints $(d-m+4,m)$ and $(\frac{d+3}{2},\frac{d+5}{2})$  
on $D_{d+1}$, with its value at $(d-m+5,m-2)$ and 
$(\frac{d-1}{2},\frac{d+7}{2})$, namely the points number
 three and two from the endpoints on $D_{d}$. 
 Thus $(d-m+5,m-2)$ is admissible only if
 $$c(d-m+5,m-2)<c(d-m+4,m)\quad {\rm and}\quad 
 c(d-m+5,m-2)<c(\frac{d+3}{2},\frac{d+5}{2}).$$
  The first inequality reduces to
 $4m-15<3d$, while the second one yields 
 $$5d^2+(68-16m)d+12m^2-100m+215>0,$$
 hence
 $$(5d-10m+23)(5d-6m+45)>-40.$$
 
 In this last inequality, when $d<2m-14$ the first factor is 
negative, so either the
  second factor is also negative, or they both have small value.  The 
second factor is negative if 
  $d<\frac {6}{5}m-9$, which is never satisfied together with 
$4m-15<3d$ for positive $m$.
  The second factor is positive with value $5d-6m+45<41$, while 
$4m-15<3d$ only if $m<32$.
  Again for values $m<32$, it is straightforward to check the 
assertion directly.

(iii) The point $(x,m)$ with $x<m$ is the endpoint with minimal $x$ 
of the diagonal $D_{d}$  
with $d=x+m-3$.   It has cardinality 
$xm - \frac{x(x+1)}{2}$, and is admissible only if the cardinality of 
the other endpoint of the diagonal $D_{x+m-2}$ is strictly bigger.  
If $x+m$ is even, then this endpoint is
$(\frac{x+m}{2},\frac{x+m+2}{2})$ and the inequality becomes
$$xm - 
\frac{x(x+1)}{2}<\frac{x+m}{2}\frac{x+m+2}{2}-
\frac{1}{2}\frac{x+m}{2}\frac{x+m+2}{2}=\frac{1}{8}(x+m)(x+m+2)$$
which means that 
$$(x-m+3)(5x-m-5)+15>0.$$
Now the first factor is negative unless $m-3\leq x\leq m-1$, while 
the second factor is negative when $x<\frac{1}{5}m+1$.
With  $x=\frac{m}{5}+2$ the inequality is satisfied only if $m\leq 
10$.  
Likewise, with $x=m-4$, the inequality is satisfied only if $m<10$.

If $x+m$ is odd, then the other endpoint of $D_{x+m-2}$ is
$(\frac{x+m-1}{2},\frac{x+m+3}{2})$ and the inequality becomes
$$xm - 
\frac{x(x+1)}{2}<\frac{x+m-1}{2}\frac{x+m+3}{2}-
\frac{1}{2}\frac{x+m-1}{2}\frac{x+m+1}{2}=\frac{1}{8}(x+m-1)(x+m+5)$$
which means
$$(x-m+3)(5x-m-7)+16>0.$$

The first factor is negative unless $m-3\leq x\leq m-1$, while 
the second factor is negative when $x<\frac{1}{5}(m+7)$.
Since $x+m$ is odd, $x<m-3$ means $x\leq m-5$.  If we set $x=m-5$, 
the inequality is satisfied only if $m<10$.
Likewise, if  $x=\frac{m}{5}+2$ the inequality is satisfied only If 
$m<\frac{153}{12}$, so the result follows for $m\geq 13$.
A special check reveals that the inequality does not hold for $m=11$ 
and 
$12$, either.
By concavity of the function $(x-m+3)(5x-m-7)+16$ then gives
that it is negative for all $m$ in the range  
$\frac{x}{5}+2<x<m-3$.

(iv) In this case the computation is similar. The point $(x,m-1)$ in
the next to upper-row has cardinality $c((x,m-1))=x(m-1) - 
\frac{x(x+1)}{2}$ and is
admissible only if it has lower cardinality
than the lower endpoint of the diagonal above. 
If $x+m$ is odd, then the lower end of the diagonal above is
$(\frac{x+m-1}{2},\frac{x+m+1}{2})$,
and its cardinality is $\frac{(x+m+1)(x+m-1)}{8}$. Now the condition
$$x(m-1) - \frac{x(x+1)}{2} < \frac{(x+m+1)(x+m-1)}{8}$$ 
translates to:
$$(x-m+3)(5x-m-3)+8 > 0.$$
The first factor is negative unless $m-3\leq x\leq m-1$, while 
the second factor is negative when $x<\frac{1}{5}(m+3)$.
Since $x+m$ is odd, $x<m-3$ means $x\leq m-5$.  If we set $x=m-5$, 
the inequality is satisfied only if $m<8$.
Likewise, if  $x=\frac{m}{5}+1$ the inequality is satisfied only If 
$m<10$.  By the concavity argument the result follows.

If $x+m$ is even, then the lower end of the diagonal above is
$(\frac{x+m-2}{2},\frac{x+m+2}{2})$,
and its cardinality is $\frac{(x+m-2)(x+m+4)}{8}$.
The necessary condition for $(x,m-1)$ being admissible is 
 $$x(m-1) - \frac{x(x+1)}{2} < \frac{(x+m-2)(x+m+4)}{8}.$$
This becomes $$(x-m+4)(5x-m-6)+16 > 0.$$
We insert $m-6$ which is the largest $x$-value smaller than $m-4$
making
$x+m$ even and obtain $m\leq 10$.
Likewise, we insert $x=\frac{m}{5}+2$ and obtain $m \leq 12$.
Hence the statement holds for $m \geq 13$.
A special check reveals that it holds for $m=11,12$ also.
\end{proof}

We return to the proof of Proposition \ref{leftright} and assume that 
$U$ is a Schubert 
union with spanning dimension $K$, and 
that $U$ has the maximal number of points among such unions, i.e. 
$g_{U}$ is maximal 
in the lexicographical order.
Therefore the grid $G_{U}$ contains an admissible point 
$\alpha=(x,y)$ in the 
$d(K)$-diagonal, i.e. $x+y-3=d=d(K)$. 
By Lemma \ref{admiss} it suffices to study the
following eight cases.

(a) $d\leq m-3$

(b)  $m>10$ and $\alpha=(d-m+3,m)$ with and $2\leq d-m+3 \le 
\frac{m}{5}+2$, i.e. $m-1\leq d \le \frac{6m}{5}-1$.

(c) $m>10$ and $\alpha=(d-m+4,m-1)$, with $2\leq d-m+4 \le 
\frac{m}{5}+2$, i.e. $m-2\leq d \le \frac{6m}{5}-2$.

(d) $d$ is even and $\alpha= (\frac{d+2}{2},\frac{d+4}{2})$.

(e) $d$ is odd and $\alpha= (\frac{d+1}{2},\frac{d+5}{2})$.

(f) $d$ is even and $\alpha= (\frac{d}{2},\frac{d+6}{2})$.

(g) $m=11$, and $G_U$ intersects the $d(K)$-diagonal in
$(4,9)$.

(h) $m\leq 10$.

In each case we consider the residual grid $\Delta=G_{U}\setminus 
G_{\alpha}$, and find the diagonal $D_{d^{\prime}}$ with largest 
$d^{\prime}$ that $\Delta$ intersects.  By the lexicographical 
ordering 
of $g_{U}$, the value of $d^{\prime}$ is determined in a similar 
fashion as $d=d(K)$ by the cardinality of $\Delta$ and the shape of 
the grid $G_{G(2,m)}\setminus G_{\alpha}$.
Notice that $U$ is a finite union of irreducible components that are 
all Schubert cycles.  Furthermore, the point $\alpha$ correspond to a 
Schubert 
cycle component $S_{\alpha}$ of maximal Krull dimension in $U$, and 
that a point $\beta\in\Delta\cap D_{d^{\prime}}$ correspond 
to a Schubert cycle  $S_{\beta}$ of maximal Krull-dimension among the 
rest of the irreducible components of $U$.
First of all the cardinality of $\Delta$ is $e=K-c(\alpha)$, and by 
definition of $d=d(K)$, the Krull dimension $d^{\prime}$ is at most 
$d(K)$.
%Since $g_{U}$ is maximal, $\Delta$ reaches the maximal diagonal $D_{d(K,\alpha)}$ 
%among all diagonals reached by residuals $\Delta^{\prime}=G_{U^{\prime}\setminus G_{\alpha}$, where 
%$U^{\prime}$ is a Schubert union of spanning dimension $K$ that contains the Schubert cycle $S_{\alpha}$.
%Notice that $d(K,\alpha)$ by definition is the maximal Krull dimension among the Schubert cycle components of 
%$U\setminus S_{\alpha}$.

Starting with (a), when $d\leq m-3$, then $\alpha=(1,d+2)$.  If 
$K\geq d+2$, then $d(K)>d+1=c(\alpha)$, 
contrary to the assumption, so $K=c(\alpha)=d+1$.
In particular $e=K-c(\alpha)=0$ and $U=S_{\alpha}$.

In case (b),  the grid $G_{G(2,m)}\setminus G_{\alpha}$ is 
similar to the original grid $G_{G(2,m)}$, but with $\{(x,y)| 
d-m+4\leq x<y\leq m\}$.  
Since $d(K)=d$, the cardinality $e=K-c(\alpha)$ of $\Delta$ is the 
cardinality of the first column of $G_{G(2,m)}\setminus G_{\alpha}$, 
i.e. at most $m-(d-m+5)=2m-d-5$.  
Therefore we may use the argument of (a) to 
conclude that $\Delta=\{(d-m+4,y)| d-m+4<y\leq e\}$.  Notice 
furthermore that $U$ clearly is maximal with repsect to the 
lexicographical order on $M_{U}$.

In case (c), with $\alpha=(d-m+4,m-1)$. Since $d(K)=d$, we first see 
that the cardinality 
of $\Delta$ is less than the cardinality of the upper row of 
$G_{\alpha}$, i.e.
$e=K-c(\alpha)) < d-m+4$.  Compare now the row of points $R=\{(x,m)| 
1\leq 
x<d-m+4\}$ with the column $C=\{(d-m+5,y)| d-m+5<y<m\}$, both in 
$G_{G(2,m)}\setminus G_{\alpha}$.
Notice that both have cardinality at least $e$, so that for $\Delta$ 
to 
reach the maximal diagonal $D_{d^{\prime}}$, it must be contained in 
one of these.
 The row $R$ starts on the diagonal 
$D_{m-2}$, while the columns $C$ starts on the diagonal 
$D_{2d-2m+8}$.  When $m>10$ and $d \le \frac{6m}{5}-2$, the highest 
of these diagonals is the 
first one, since then $2d-2m+8\leq \frac{2m}{5}+4 <m-2$, 
so in that case $\Delta$ must be completely contained in the row $R$.

To see whether $U$ is maximal with respect to the lexicographical 
order 
of $M_{U}$, there are essentially two different situations:
$e=d-m+3$ (maximum possible), and $e \le d-m+2$. If $e=d-m+3$, we are
already in case a), since $(d-m+3,m)$ and $(d-m+4,m-1)$ are on the 
same
diagonal.

Assume $e \le d-m+2$. Then we collectively remove the $d-m+3-e$ top 
squares
of the right column of $G_{\alpha}$, and reinstall them horizontally
as points $(e+1,m), (e+2,m)...,(d-m+3,m)$. This amounts to moving 
squares along  $d-m+3-e$ diagonals, and does not alter the number
of $F_q$-rational points of the Schubert unions represented
by the two grids. But after moving, we have a union which is maximal
with respect to the lexicographic order on the $M_U$, and we are 
done. 
(We have $d-m+3$ columns to the left filled up completely,
and $K-c(d-m+3,m)$ squares in column nr $d-m+4$).

Starting with (d), with $\alpha= (\frac{d+2}{2},\frac{d+4}{2})$  
the first row in $G_{G(2,m)}\setminus G_{\alpha}$ has length larger 
than $d^{\prime}$, since $d(K)>d^{\prime}$.  Therefore $\Delta$ is 
completely contained in this first row.  Furthermore $U$ is clearly 
minimal with respect to the lexicographic order of $M_U$.

Starting with (e) we see that $K-c(\alpha)=0$, since if $\Delta$ is 
non-empty, then it could lie to the right of $\alpha$ in the top row, 
and $d(K)$ would have been larger. On the other hand $U=S_{\alpha}$ 
is 
clearly minimal with respect to the lexicographic order of $M_U$.

Starting with (f), we see that $K-C(d(K))=0$ or $1$, since if it was
at least $2$, $\Delta$ could lie to the right of $\alpha$, and
then $d(K)$ would have been larger.
If $\Delta$ is empty, there is nothing to do, and if it is one point 
the optimal solution is that it lies to the right of $\alpha$.
In both cases $U$ is minimal
with respect to the lexicographic order of $M_U$.

Starting with (g), we have $m=11$, and $\alpha=(4,9)$. Since $\alpha$
lies on the diagonal $D_{10}$ and
$c(4,9)+1=27=C(11)$, we see that this is only an issue when 
$K=c(4,9)=26$, i.e. $Delta$ is empty. But a quick calculation reveals 
that
$S_(4,9)$ is not optimal at all for $K=26$. An optimal choice in this
case is  $U=S_{(2,11)} \cup S_{(3,10)}$, with $M_U=\{10,9,7\}$, which
is maximal with respect to the lexicographical order.
Hence the assertion holds in all cases.

The cases $m\leq 10$ are entirely similar. The check that no case 
occurs that is not of one of the kinds above is left to the reader.
\end{proof}

\begin{rem} \label{lr}
{\rm From Proposition \ref{leftright} it is clear that for each 
spanning 
(co)dimension we only need to check two Schubert unions to find one 
which is maximal with respect to $g_U$. In the tables below we 
utilize this fact,
and indicate with an $L$ (go left) if we may use the union maximal 
with respect 
to $M_U$, with an $R$ (go right) if we may use the union minimal with 
respect 
to $M_U$, and with $LR$ if and only if we may use both.
The spanning codimension is $r=$$m \choose\ 2$$-K$.}
\end{rem}

C(2,7): 
\[ \left[
\begin{array}{cccccccccccc} 
{\rm Codim:}  &0       & 1      & 2   & 3   & 4   & 5   & 6& 7   & 
8   & 9   & 10   \\
{\rm Direction:} &LR & LR & LR&LR & R & R &R &R  & R & LR& L 
\end{array}
\right]\]
\[ \left[
\begin{array}{cccccccccccc} 
{\rm Codim.:}   & 11  & 12   & 13 & 14  & 15  & 16  & 17  & 18  & 19  
& 20 & 21   \\
{\rm Direction:} &R & LR  & L & L & L&L &L &LR &LR &LR &LR
\end{array}
\right]\]

%C(2,8): 
%\[ \left[
%\begin{array}{cccccccccccccccc} Codim.: &0& 1 & 2 & 3 & 4& 5 &  6 & 7 & 8 & 9 & 10& 11 &12&13 &14 \\
%Direction: &LR& LR & LR &LR & R & R &R & R  &R & R & LR &R&R&R&LR
%\end{array}
%\right]\]

%\[\left[
%\begin{array}{ccccccccccccccc} 
%Codim.: & 15  & 16  & 17  & 18  & 19 &  20 & 21 & 22  & 23  &
%24  & 25  & 26  & 27 & 28 \\
%Direction:& L &L & L & LR & L & L  &L &L&L &L &LR& LR& LR &LR
%\end{array}
%\right]\]

C(2,9): 
\[ \left[
\begin{array}{cccccccccccccc} {\rm Codim.:} &0& 1 & 2 & 3 & 4& 5 &  6 
& 7 & 8 & 9 & 10& 11 &12 \\
{\rm Direction:} &LR& LR & LR &LR & R & R &R & R  &R & R & R &R&R
\end{array}
\right]\]

\[\left[
\begin{array}{ccccccccccccc} 
{\rm Codim.:} &13&14& 15  & 16  & 17  & 18  & 19 &  20 & 21 & 22  & 
23  &
24  \\
{\rm Direction:} & R &R & R & R & R & LR  &L &L&L &L &L& L
\end{array}
\right]\]

\[\left[
\begin{array}{ccccccccccccc} 
{\rm Codim.:} &25&26& 27  & 28  & 29  & 30  & 31 &  32 & 33 & 34  & 
35  &
36  \\
{\rm Direction:}& L &L & L & L & L & L  &L &L&LR &LR &LR& LR
\end{array}
\right]\]

C(2,10): 
\[ \left[
\begin{array}{ccccccccccccccccc} {\rm Codim.:} &0& 1 & 2 & 3 & 4& 5 
&  6 & 7 & 8 & 9 & 10& 11 &12&13&14&15 \\
{\rm Direction:} &LR& LR & LR &LR & R & R &R & R  &R & R & R 
&R&R&R&R&R
\end{array}
\right]\]

\[\left[
\begin{array}{ccccccccccccccc} 
{\rm Codim.:} & 16  & 17  & 18  & 19 &  20 & 21 & 22& 
23&24&25&26&27&28&29 \\
{\rm Direction:} & R &R & R & R & LR & L  &R &R&R &LR &L& L&L&L
\end{array}
\right]\]

\[\left[
\begin{array}{ccccccccccccccccc} 
{\rm Codim.:} & 30& 31 &  32 & 33 & 34  & 35  &36 
&37&38&39&40&41&42&43&44&45 \\
{\rm Direction:} &L& L &L & L & L & L & L  &L &L&L &L &L& LR&LR&LR&LR
\end{array}
\right]\]

Each table starts and ends with $4$ occurences of $LR$.
This is because in the three largest and the three smallest spanning
dimensions
$K$ there is only one Schubert union, and because we have only two 
Schubert unions with spanning dimension $3$, namely  $S_{(2,3)}$, in 
projective terms
a $\beta$-plane, or $S_{(1,3)}$, an  $\alpha$-plane. Both have
$q^2+q+1$ points. In codimension $3$ we have the duals of these two,
of course also with the same number of points. 

From the tables for $G(2,8)$ (which is not listed) and $G(2,9)$ one 
can
conclude without further computations that the $E_r$ are always 
monomials of type $q^i$ in these cases (See Question (Q7) of
Section \ref{w67}). That is because we never
jump directly from an $R$ to an $L$ or vice versa in these cases, we
always go via an $LR$. For $m=7$ there is a jump between 
$L$ and $R$ between codimensions $10$ and $11$, but a calculation
reveals that $J_{10}-J_{11}=q^5$.  For $m=10$ we observe the fatal 
jump from $R$ to $L$, passing from codimension $22$ to $21$. 
Here $E_{22}$ is not a monomial in $q$. 
As opposed to the tables above it, the one for $(l,m)=(2,10)$ is not 
symmetric in $L$ and $R$.

\begin{cor} \label{twocycles}
For $G(2,m)$ and any $m$: There are no Schubert unions that are proper
unions of more than two Schubert cycles, and that contain the maximal
number of $F_q$-rational points, given its spanning dimension.
\end{cor}
\begin{proof}
This is not a corollary of Proposition \ref{leftright}, which only
says that for each spanning dimension $K$, there exists a Schubert
union of a special kind  containing the maximum number of points.
On the other hand the proof of Proposition \ref{leftright} gives the
corollary as follows: For an optimal Schubert union we have to be in 
one of the situations (a),(b),(c),(d),(e), (f) or (g) described in 
that proof. In each case the refined analysis gave
that an optimal union would have to look in a specific way.
Tracing the arguments, we see that it is obvious that the optimal
unions
are proper unions of one or two cycles in all cases (a), (b), (c), 
(d), (e) and 
(f).

Case (g) ($m=11, K=26$) is checked explicitly.
\end{proof}

\begin{rem} \label{sumup}
{\rm 
One might ask: For a given  $l,m$ and spanning dimension $K$:
Do there exist Schubert unions with a maximal number of points
amomg those of that spanning dimension that are proper Schubert unions
of at least $l+1$ Schubert cycles. Corollory \ref{twocycles} says:
``No, if $l=2$''.
On the other hand it is clear that if $l \ge 3$, then there are 
optimal
Schubert unions that are proper unions of more than two Schubert
cycles.
We will see, in the tables for $G(3,6)$ in Section \ref{Tables} that
$S_{(1,5,6)} \cup S_{(2,3,6)}\cup S_{(3,4,5)}$,
is among the ones with a maximal number of points for $G(3,6)$ and 
$K=15.$}
\end{rem}

We now will investigate, for each $m$, the range of those spanning
dimensions $K$, where we will use the Schubert union $U$ with maximal
$M_U$ to compute the maximum number of $F_q$-rational points, and the
range
where we will use the $U$ with minimal $M_U$ to do the same.
We will give some general results. In fact, we are well on our way,
through Lemma \ref{admiss} above. We never used part (iii) of that 
result to
prove Proposition \ref{leftright}, but we will use part (iii) now.

\begin{prop} \label{whenleft}
(i) If $d(K) > 1.2m - 1$, then the Schubert union $U$ with minimal
$M_U$ will have a maximal number of $F_q$-rational points.

(ii)  If $d(K) \le 1.2m - 5$, then the Schubert union $U$ with maximal
$M_U$ will have a maximal number of $F_q$-rational points. 
\end{prop}
\begin{proof}
In case (i) we will prove that the union with maximal $M_U$ is
impossible, unless $d(K) \ge 2m-7$, and for these $d(K)$ the matter
is obvious. In case (ii) we will prove that the union with minimal
$M_U$ is impossible, for diagonals that intersect the top row of
$G_U$. For the remaining ones the result is obvious.

In case (i), to be able to use the union with
maximal $M_U$, a necessary condition is that the point $\alpha=(x,m)$ 
where
$G_U$ intersects the diagonal $x+m-3=d(K)$ is admissible.
But from  Lemma \ref{admiss}, (iii) we know that $\alpha$ is
admissable only 
if $x \le 0.2m+2$, or $x \ge m-3$. Hence, amomg the two unions in
question,  only the $U$ with
minimal $M_U$ is possible for 
$d(K) > (0.2m+2) +(m-3)=1.2m-1.$ and $d(K)\le (m-4) +m-3=2m-7$.
But for $d(K) \ge 2m-6$ a quick glance at a diagram reveals that the
$U$ with minimal $M_U$ is optimal. 

In case (ii) we perform the same sort of calculation. To handle
diagonals
 that intersect the top row we compare
the cardinality $(x+1)m-\frac{(x+1)(x+2)}{2}$ of a point $(x+1,m)$ in 
the top
row with that of the bottom-right point in the diagonal below it.
That cardinality is $\frac{(x+m+1)(x+m-1)}{8}$ if we have a point 
which is 
rightmost in its row, and  $\frac{(x+m+2)(x+m)}{8}-1$  if it isn't.
In the first case we obtain that the cardinality $c(x+1,m)$ is the 
smaller one if
$x \le 0.2m-1.4$ which gives $d(K)=(x+m-3 \le 1.2m -4.4$.
In the second case we obtain that $x \le 0.2m-1.4$, corresponding to 
 $d(K)=(x+m-3 \le 1.2m -5$ is enough
For diagonals that intersect the left side of $G_{G(2,m)}$ the result
is obvious.
\end{proof}
\begin{rem} \label{estimates}
{\rm (i) The results in Proposition \ref{whenleft}
are not best possible, direct calculations for
each $m$ will often give results that are a little sharper. 
The rough, but essential, picture is that we change from Schubert 
unions that are maximal 
to unions that are minimal, with respect to $M_U$, 
when we reach a diagonal which intersects the upper edge of the
$G_{G(2,m)}$ about $20\%$ of the way from the left corner to the right
one.

(ii) If one prefers bounds on $K$ instead of on $d(K)$, one can do as 
follows. The argument above showed: For $x \le 0.2m - 2$ 
(corresponding to
$d \le 1.2m-5$) the points 
$(x,m)$ in the top row of $G_{G(2,m)}$, are contained in $G_U$ for
  unions $U$ which are optimal for spanning dimensions $K$ with
$d(K)=x+m-3$. 
%The points at the other ends are not contained in such $G_U$.
Let $a$ be the integral value $[0.2m-2]$. Then we can fill up the
$a$ first columns with squares, and be certain 
that the grids thus formed, are associated with unions, maximal with 
respect to $g_U$. This gives:
For $K \le(m-1)+(m-2)+....+(m-a)=am - $$a+1 \choose 2$
we can always find maximal unions with respect to $g_U$ that are 
maximal with respect to $M_U$. 
Since $0.2m-3 \le a \le 0.2m-2$, this number is at least
$0.18m^2-2.7m-1$. But $k= $$m\choose2$ is $0.5m^2-0.5m,$ so we see
that for about $36\%$ of the $K$ we can be sure to ``stick to the 
left''
when $m$ is big.
But a similar argument can be used, using Proposition \ref{whenleft}, 
(ii),
and then we see that for about $64\%$ of the $K$ we can be sure to
``stick to the right''.}
\end{rem}
This gives rise to the following result for the $G(2,m)$:
\begin{prop} \label{percentage}
For every $\epsilon > 0$, there exists a natural number  $M$, such 
that
if $m >M$, then 

(i) If  $K \le 0.36k-\epsilon$, then he $U$ with maximal $M_U$
is maximal with respect to $g_U$, and  
the $U$ minimal with respect to $M_U$ is not maximal with respect to 
$g_U$ 
unless $K \ge k-3$.

(ii) If  $K \ge 0.36k+\epsilon$, then he $U$ minimal with respect to 
$M_U$
is maximal with respect to $g_U$,
and  the $U$ maximal with respect to  $M_U$ is not maximal with 
respect to 
$g_U$unless $K \le 3$. 

\end{prop}

\begin{rem} \label{nondual}
{\rm  (i) Let $l=2$. For every $K$ between $0$ and $k= $$m \choose\ 
2$ we know that 
if $U$ has spanning dimension $K$ and is  minimal with respect to 
$M_U$, then the
dual of $U$ has spanning dimension $k-K$ and $U^*$ is maximal with 
respect to
$M_{U}$. From Proposition \ref{percentage} we see that in the range 
between $0.36k+\epsilon$ and
$0.64m-\epsilon$, say in the interval $[0.37k,0.63k]$, the $U$ which 
is minimal
with respect to $M_U$, but not its dual, will be maximal with respect 
to 
$g_U$ if $m$ is large enough. 
For $l=2$, the smallest $m$ such that there exists an $K$ with a 
Schubert union
$U$, minimal or maximal with respect to $M_U$ and maximal with 
respect to $g_U$, such that
its dual $U^*$ is not maximal with respect to $g_U$ (for spanning 
dimension $k-K$), is $m=10$ (and $K=22$). 
%At the end of this section, for $l=2$ and $7 \le m \le 10$  we give some 
%explicit tables indicating whether the union maximal, or 
%the union minimal, with respect to $M_U$, is maximal with respect to $g_U$.  

(ii) Set $l=2$. Proposition \ref{percentage} gives that the percentage
of the $r$ such that $E_r$ is not of the form $q^i$, goes to zero, as
$m$ grows to infinity, ($r$ must be very close to  $0.64k$,
for smaller $r$ we ``stick to the left side'', for bigger ones we 
``stick to the right side'').}
\end{rem}

\begin{rem} \label{asymptotic}
{\rm Let us interpret  Proposition \ref{percentage} in a continous 
setting. 
Study the triangle with corners $(0,0),
(0,1),(1,1)$. This is the ``limit'' of $G_{G(2,m)}$ scaled down a
  factor $m$ in all directions as m goes to infinity.
Look at the tetragon $T_x$ with corners $(0,0),(x,x),(x,1), (0,1)$.
This tetragon has area $A=x-x^2$. We also study the triangle $P_y$ 
with
corners $(0,0),(y,y),(0,y)$ and area $A=\frac{y^2}{2}$.
Reversing the area formulas we get $x=1-\sqrt{1-2A}$ for the tetragon,
and $y=\sqrt{2A}$ for the triangle.
Now we look at diagonals of the form $x+y=d$.
The largest $d$ for which the trapes $T_x$ intersects such a diagonal
is $d_1(A)=1+x=2-\sqrt{1-2A}$, where $A$ is the area of the trapes.
The largest $d$ for which the triangle $P_y$ intersects such a 
diagonal
is $d_2(A)=2y=2\sqrt{2A}$, where $A$ is the area of the triangle.
The decisive criterion for whether it is optimal (with respect to the
lexicographical ordering on the $g_U$) is ``how much $d$'' you can
obtain with a given area $A$.
It is now easy to see that $d_1(A) > d_2(A)$ for $0 \le A < 0.18$,
and the values are equal for $A=0.18$ and $A=0.5$, while
$d_1(A) < d_2(A)$ for $0.18 < A < 0.5$. Since the area of the whole
triangle is $\frac{1}{2}$, the (exact) value $0.18$ of course
corresponds to $36\%$ of the area. }
\end{rem}

%\eject

\section{Some questions and answers about the Grassmann codes 
$C(l,m)$}
\label{w67}

For fixed $l,m,q$ let $C(l,m)$ be the Grassmann code over $F_q$ 
described in Section \ref{codes}. Recall the invariants $H_r, 
\Delta_r, J_r,D_r,E_r$ introduceed in Definitions \ref{Hr}, 
\ref{Sch}, and \ref{Sch2}.
inspired by Proposition \ref{Nogin}, Theorem \ref{highweights}, 
Example
\ref{CIproj}, Proposition \ref{d5}, and the results of Section
\ref{algo} we now will formulate some natural questions, which we 
will also comment on briefly:

%\subsection{Questions and observations}
%\label{questions}

For each $l,m,q$  we obviously have: 
\begin{equation} \label{simply}
\sum^{k}_{r=1}\Delta_{r}=d_k=
n.
\end{equation}
Here $n$ and $k$ are the word length and dimension of $C(l,m)$ as 
before. 
Moreover it is clear that $n$ is the sum of $k= $$m \choose\ l$
monomials of type $q^i$.
For each $l,m$ one may raise the following questions:

\begin{itemize} \label{questions}

\item (Q1) Are the $d_r$  always sums of $r$ monomials of type
$q^i$, for $r=1,...,k$ ?

\item (Q2) Is  $\Delta_{r}$ always a monomial of the form $q^i$ ?

\item (Q3) Is it true that:
%To our knowledge, all known results indicate that this holds, for all
%$C(l,m)$ and all $r$. We conjecture that (*)  always holds, and
%believe that (**) also always holds. 

%In addition, all known results for $l=2$ indicate that we have

$$ \Delta_r(q)=q^{l(m-l)}\Delta_{k+1-r}(q^{-1}),$$
for all $C(l,m)$, and all $r$ ?
This in turn implies that if  the answer to question (Q2) is (partly)
positive,
and $\Delta_{r}=q^i$ for some $i$, then 
$\Delta_{k+1-r}=q^{l(m-l)-i}$.

\item Answers to (Q1), (Q2), (Q3): Affirmative for 
$(l,m)=(2,3),(2,4),(2,5)$
by Proposition \ref{Nogin}, Theorem \ref{highweights}, and 
Proposition \ref{d5}. In other cases we do not know the answers for 
all $r$
(affirmative for the smallest and biggest $r$).

Question (Q3) is inspired by Proposition
\ref{dualnumbers} and the following question:

\item (Q4) If $H_r$ is computed by a linear section $H$ of $G(l,m)$, 
is then 
$H_{k-r}$ always computed by $D(H)$ in the sense of Definition 
\ref{duality}?  

\item Answer: We do not know.

\item (Q5) Is it true that $J_r=H_r$, and therefore $D_r=d_r$, and 
$E_r=\Delta_r$,
for all $l,m,r$ ? 

\item Answer: Affirmative for $(l,m)=(2,3),(2,4),(2,5)$.
In other cases we do not know the answers for all $r$
(affirmative for the smallest and biggest $r$).

Taking into account the possibility that the answer to question
(Q5) is no, we may phrase similar questions as (Q,1-4) 
with the $J_r, D_r, E_r$ replacing $H_r, d_r, \Delta_r$, respectively:

\item (Q6) Are the $D_r$ and $J_{r}$ always sums of $r$ monomials of 
type
$q^i$, for $r=1,...,k$ ?

Answer: Affirmative for all $(l,m)$ by Proposition \ref{int-un}.

\item (Q7) Is  $E_{r}$ always a monomial of the form $q^i$ ?

\item Answer: Affirmative for $(l,m)=(2,m)$, for $m \le 9$, and 
$(l,m)=(3,6)$.
Negative for some $r$ for $(2,m)$, and $m=10$ or $m$ big enough.
We have not performed further investigations.

\item (Q8) If $J_r$ is computed by a Schubert union $S_U$, is 
$J_{k-r}$ then computed by $S_{U^*}$? 

\item Answer: Affirmative for $(l,m)=(2,m)$, for $m \le 8$, and 
$(l,m)=(3,6)$.
Negative for $(2,m)$, and $m=10$ or big enough.
We have not performed further investigations.

\item (Q9) Is it true that:

$$E_r(q)=q^{2m-4}E_{k+1-r}(q^{-1}) ?$$

for all $C(l,m)$, and all $r$ ?
This in turn would imply that if the answer to question (Q7) is 
(partly)
positive, and $\Delta_{r}=q^i$ for some $i$, then 
$\Delta_{k+1-r}=q^{l(m-l)-i}$.

\item Answer: Affirmative for $(l,m)=(2,m)$, for $m \le 9$, and 
$(l,m)=(3,6)$.
Negative for $(2,m)$, and $m=10$.
We have not performed further investigations.
\end{itemize}

\begin{rem} \label{sum}
{\rm It follows from the results of Section \ref{codes}
that all questions, except possibly (Q3) have affirmative answers for 
$l=2$ and $m \le 5$. The affirmative answers to  (Q7), (Q8), (Q9) for 
$(l,m)=(2,6),(2,7),(2,8),(3,6)$ are due to explicit investigations. 
See the tables in the appendix.  For $(l,m)=(2,9)$ it
is at least clear that (Q7) and (Q9) have affirmative answers. See 
Remark \ref{lr}. The
negative parts of the answers to these questions follow essentially 
from 
Remark \ref{nondual}. See Remark \ref{nondual}(ii) though, concerning 
(Q7).
For $(l,m)=(2,10)$ we see from Remark \ref{lr} and explicit 
calculations that  $E_{22}=J_{21}-J_{22}=q^9+q^8-q^6$,
so (Q7) has a negative answer. Moreover $E_{24}=J_{23}-J_{24}=q^6$, 
and hence
(Q9) also has a negative answer.

Detailed descriptions of Schubert unions for
$(l,m)=(2,4),(2,5),(2,6),(3,6)$ are given in the Appendix as
illustrations. }
\end{rem}

From the observations above we may conclude: 

\begin{prop} \label{No}
Neither of the questions question (Q7), (Q8), and (Q9) do always have 
affirmative answers, and questions
(Q1), (Q2), (Q3), (Q4), and (Q5) do therefore not 
simultaneously have affirmative answers for all $l,m,r,q$.
\end{prop}

\section{Linear section with maximal Krull dimension}
\label{Krull}

Inspired by the questions in the preceeding paragraphs we may ask $3$
addition questions for fixed $(l,m)$, and $K$ in the range between 
$1$ and $k$:
\begin{itemize} \label{quest}

\item (q1) What is the maximal Krull dimension of a component of a 
linear section of $G(l,m)$ with a subspace of Pl\"ucker space of 
(projective) dimension 
$K-1$ ?

\item (q2)  What is the maximal Krull dimension $d(K)$ of a component 
of a linear section of $G(l,m)$ with a subspace of Pl\"ucker space of 
(projective) dimension $K-1$, spanned by a Schubert union ?

\item (q3) Are the answers to (q1) and (q2) identical ?

\end{itemize} 

The last question is obviosly an analogue of,  and weaker version of, 
(Q5) from
Section \ref{w67}. It is clear that the analysis in Section 
\ref{algo}, in particular the proof of Proposition \ref{leftright}, 
gives the answer to (q2) for $l=2$.

For each $d \ge m-2$ let $c_1(d)$ be the cardinality of the upper 
point of the
diagonal corresponding to Krull dimension $d$, that is of 
$(d-m+3,m)$. This is: 
$$c_1(d)=(d-m+3)m-\frac{(d-m+3)(d-m+4)}{2}=\frac{4dm-3m^2-d^2+13m-7d-12}{2}.$$

Moreover we let $c_2(d)$ be the cardinality of the lower point of the
diagonal corresponding to Krull dimension $d$. This point is 
$(\frac{d+2}{2},\frac{d+4}{2})$ if $d$ is even and 
$(\frac{d+1}{2},\frac{d+5}{2})$ if $d$ is odd. Hence 

$c_2(d)=\frac{d^2+6d+8}{8}$, if $d$ is even, and 
$c_2(d)=\frac{d^2+8d+7}{8}.$ if $d$ is odd.
Furthermore we set 

$$C(d)=\min\{c_1(d),c_2(d)\}.$$ For $d \le m-2$ we set $C(d)=d+1$.

We obtain from the arguments in Section \ref{algo}: 
\begin{prop} \label{Schubert-Krull}
 The maximal Krull dimension $d(K)$ of a component of a linear 
section of $G(l,m)$ with a subspace of Pl\"ucker space of 
(projective) dimension 
dimension $K-1$, spanned by a Schubert union, is the largest $d$ such 
that $C(d) \le K$. 
\end{prop}
Now we will argue that the answer to question (q3) is affirmative, 
thereby also answering question (q1):
Recall the notation from Propositions \ref{Borel1} and \ref{Borel2} 
and the text between these two results.
For $1 \le r \le k$ let $G(r,V_l)$ be the Grassmann variety 
parametrizing 
projective $(r-1)$-spaces in the Pl\"ucker space $\bf P$ containing 
$G(l,m)$. Let $I$ be the incidence variety in $G(l,m) \times 
G(r,V_l)$ parametrizing inclusion relations of points in $G(l,m)$ and 
linear sections of $\bf P$ corresponding to points of $G(r,V_l)$. Let 
$f$ be the projection to the second factor.
Let $E$ be the subset of $I$ corresponding to those $x$ of $I$ 
whose inverse image $f^{-1}(f(x))$ contains a component of maximal 
dimension 
among the fibres of $f$. By Exercise II, 3.22 of \cite{H}, we see 
that $E$ is 
closed in $I$. Moreover the set $F=f(E)$ is closed in $G(r,V_l)$ by 
\cite{H}, 
Exercise II, 4.4, since $E$ is closed and proper. By the argument 
preceeding 
Proposition \ref{Borel2}
the solvable Borel group $B$ acts on $V_{l,r}$, and the set $F$ is 
closed in 
$P(V_{l,r})$ in virtue of being closed in  $G(r,V_l)$.
The set $F$ is stable under the action of 
$B$ since the property of having a component with maximal fibre 
dimension under $f$
is invariant under the action of $B$.
Since $B$ is irreducible each irreducible component of $F$ is 
therefore
stable under $B$. Hence $B$ acts on each component of $F$, which is 
closed in   $P(V_{l,r})$ where the action is induced by a linear 
action on  $V_{l,r}$.
Borel's Fixed Point Theorem, as quoted for eample in \cite{FH}, p. 384
 (or Theorem 7.2.5 of \cite{Sp}, or \cite{Fu}, p. 155, see also Remark \ref{closed}
 below),
then gives that $B$ must have a fixed point in each component of $F$, 
so at least one fixed point. But by Proposition \ref{Borel2} the 
fixed points correspond precisely to the projective $(K-1)$-planes 
spanned by Schubert unions.

Hence at least one of these special $(K-1)$-planes have a component 
of maximal fibre (Krull) dimension under $f$. This gives an 
affirmative answer to question [q3), and also gives the following 
result for $l=2$:

\begin{thm} \label{maxkrull}
 The maximal Krull dimension $d(K)$ of a component of a linear 
section of $G(l,m)$ with a subspace of Pl\"ucker space of 
(projective) dimension 
dimension $K-1$, is the largest $d$ such that $C(d) \le K$. 
\end{thm}
\begin{rem} \label{closed}
{\rm  Borel's Fixed Point Theorem is valid for algebraically closed fields $F$, 
and we have now shown that for such fields the linear sections of $G(l,m)$ with a
component of the maximal possible Krull dimension can be found among Schubert unions. These are given by linear equations of type $X_{i,j}=0$, in other words defined over ${\bf Z}$ and hence over any subfield of $F$ (here we take only those linear sections given by equations with coefficients in the subfield). Hence Theorem \ref{maxkrull} is valid also for non-algebraically closed fields.
}
\end{rem}
The affirmative answer to question (q3) is an indication that the 
following conjecture holds.
\begin{conj} \label{Schubertmax} 
The answer to Question (Q5) of Section \ref{w67} is always 
affirmative.
Hence the higher weights of the Grassmann codes are always computed 
by Schubert unions.
\end{conj}

\section{Codes from Schubert unions}
\label{unionskoder}

In earlier sections we have studied the impact of Schubert unions
to Grassmann codes in order to make the bound $d_r \le D_r$ explicit.
Now we will study codes made from a Schubert union $S_U$ in the same
way as the codes $C(l,m)$ are made from the $G(l,m)$. In other words;
For a given Schubert union $S_U$ and prime power $q$ 
denote the (affine) spanning dimension of $S_U$ by $K_U=K$.
Then the Pl\"ucker coordinates of all points of $S_U$ have only zeroes
in all the coordinates corresponding to the $k-K$ points of $H_U$,
so we delete them. Choose coordinates for each points, and make the
corresponding $K$-tuples columns of a $k \times g_U(q)$-matrix $\G$.
This matrix will be the generator matrix of a code.
If we change coordinates for a point by multiplying by a factor, 
the code changes, but its equivalence class and code parameters do
not, so by abuse of notation we denote  all equivalent codes appearing
this way bu $C_U$.

In \cite{HC} it was shown that if $l=2$, and we simply have a Schubert
cycle $S_{\alpha}$, then the minimum distance $d_1=d$ of the code is
$q^{\delta}$, where $\delta$ is the Krull dimension of the
Schubert cycle. We will use this result to give the following
generalization:
\begin{prop} \label{mindist}
For a Schubert union $S_U$ in $G(2,m)$, which is the proper union
of $s$ Schubert cycles $S_i$ with Krull dimensions $\delta_i$, for
$i=1,..,s$, the minimum distance of 
$C_U$ is the smallest number among the $\delta_i$.
\end{prop} 
\begin{proof}
Let $S_{\alpha}$ be one of the cycles in the given union with minimal 
Krull dimension $\delta$. We now intersect $S_U$ with the coordinate
hyperplane $X_{\alpha}$ (restricted to the $K$-space in which 
$S_U$ sits, if one prefers). Since $\alpha$ is not contained in the
$G_{\beta}$ of any Schubert union $S_{\beta}$
different from $S_{\alpha}$ appearing in the union, this coordinate
hyperplane contains all these $S_{\beta}$. By standard arguments
there are exactly $q^{\delta}$ points from $S_{\alpha}$ that are not
contained in this hyperplane (all these points are then of course 
outside all
the other $S_\beta$). Standard argument: If $\alpha=(a,b)$, then this
hyperplane cuts out $S_{(a-1,b)} \cup S_{(a,b-1)}$, with exactly one 
point $(a,b)$ less in its $G$-grid, and by Corollary
\ref{int-un} we must then subtract $q^{a+b-3}$ to obtain the number
of points.
On the other hand it is clear that if we intersect $S_U$ with an
arbitrary
hyperplane $H$ in $K$-space (or an arbitrary hyperplane in the 
Pl\"ucker space,
not containing $S_U$), then there is at least one $S_i$, which is not
contained in $H$. Now the maximal number of points of any hyperplane
section of $S_i$ is equal to the cardinality of $S_i$ minus
$\delta_i$, so there are at least $q^{\delta_i}$ points of $S_i - H$.
Hence there are at least $\delta_i$ points of $S_U - H$ also.
Hence the maximal number of points of $S_U \cap H$ is $g_U(q) -
\delta$,
where $\delta$ is the smallest $\delta_i$, and $d=d_1$ is computed by
$X_{\alpha_i}$ for such a corresponding $i$. 
\end{proof}

We may also  mimick the contents of 
Proposition \ref{Nogin} and Theorem \ref{highweights}. 
Let $\alpha$ be such that $S_{\alpha}$ is one of the Schubert cycles
$S_i$ with minimal Krull dimension in $S_U$, and set $\delta=\delta_i$
(not necessarily the degree of $g_U$).
\begin{prop} \label{Nogin2}
$d_r=q^{\delta}+q^{\delta -1}+ \dots +q^{\delta -r+1}$, for 
$r= 1, \dots s$, where $s$ is largest  natural number such that 
$(a-s+1), (a-s+2,b),...,(a-1,b), (a,b)$ all are contained in $G_U$.
\end{prop}
\begin {proof}
$d_r \ge q^{\delta}+q^{\delta -1}+ \dots +q^{\delta -r+1}$ is the
Griesmer bound. The opposite inequalities follow if we can exhibit
linear spaces with increasing codimension, which intersect $S_U$
in an appropriate number of points. We intersect with:
  $$X_{(a,b)}=X_{(a-1,b)}= X_{(a-2,b)}=....=X_{(a-r+1,b)}=0$$
Then, as intersections  we obtain smaller successive Schubert unions. 
Their cardinalities are determined by Corollary \ref{int-un} 
and the fact that we peel off points one by one to obtain the
successive $G$-grids.
\end{proof}
We also have: 
\begin{prop} \label{highweights2}
Let $b=b_s$ be the largest number appearing in the
sequence $\sigma_U = a_1 < a_2 <.... <a_s <b_s <....<b_2<b_1.$
Then  $d_K=g_U(q)$, and
\[d_{k-a}=g_U(q)-(1+q+ \dots q^{a-1}),\]
for $a=1, \dots ,b_s-1$.
 \end{prop}
\begin{proof}
$S_U$ contains a projective space of dimension $b_s-2$.
See the proof of Theorem \ref{highweights}.
\end{proof}
Of course we also have a relative bound of type $d_r \le D_r$
\begin{prop}  \label{relSchub}
Let $S_U$ be a Schubert union in $G(l,m)$, and let $M_r$ be the
maximum cardinality of a Schubert union that is contained in $U$,
and whose spanning dimension is $r$ less than that of $S_U$.
Then $d_r \le g_U(q)-M_r$. 
\end{prop}
The proof is obvious.
\begin{exa} \label{exa}
In the apppendix (Section \ref{Tables}) we will list Schubert unions 
that
compute the $d_r$ for the Grassmann code $C(2,5)$ from $G(2,5)$.
We leave it to the reader to find the full weight hierarchy for all
$C_U$, for all $15$ non-empty Schubert unions $U$ of $G(2,5)$, using
the results above and  the table for $G(2,5)$ in the appendix.
\end{exa}

For $l\ge3$ the expected result $d=d_1=q^{\delta}$ for Schubert cycles
has not yet been shown. If it is shown, we see that we can extend it 
to Schubert unions as in the case $l=2$, and also a variant of
Proposition \ref{Nogin2} will then follow. A variant of Proposition 
\ref{highweights} holds already, and we leave it to the reader to
formulate it.

\vfill
%\eject
\section{APPENDIX: TABLES OF SCHUBERT UNIONS}
\label{Tables}
We first give  tables of the $E_r$ 
for  $C(2,m)$, for  $m=6, 7, 8$, and for $C(3,6)$.
The values are determined 
using Corollary \ref{int-un}.
% and going  through all
%$2^{m-1}$ Schubert unions for  $C(2,m)$ for low $m \le 8$,
%and all $66$ Schubert unions for $C(3,6)$. 
%For other ways to
%investigate optimal Schubert unions, see Section \ref{algo}.
%Later in this section we also give lists of all Schubert unions in
%some of these cases. 

\[C(2,6):  
\left[
\begin{array}{cccccccccccccccc} 
r:         & 1   & 2   & 3   & 4   & 5   & 6 &7&8 &  9   & 10  & 11  
& 12& 13 & 14 & 15       \\
E_r:  & {\bf q^8} & {\bf q^7} & {\bf q^6 } & {\bf q^5} & {\bf q^4} &
q^6 & q^5&q^4 & q^3 & q^2 & {\bf q^4} & {\bf q^3} & {\bf q^2} & {\bf 
q} &
        {\bf 1}  
\end{array}
\right]\]

$C(2,7):$
\[ \left[
\begin{array}{ccccccccccc} 
r:         & 1      & 2   & 3   & 4   & 5   & 6& 7   & 8   & 9   & 
10    \\
E_r:  & {\bf q^{10}} & {\bf q^9} &{\bf  q^8} & {\bf q^7} & {\bf q^6}
&{\bf q^5} &q^8  & q^7 & q^6 & q^4 
\end{array}
\right]\]
\[ \left[
\begin{array}{ccccccccccc} 
 11  & 12   & 13 & 14  & 15  & 16  & 17  & 18  & 19  & 20 & 21   \\
q^5 & q^6  & q^4 & q^3 & q^2 &{\bf
  q^5} &{\bf q^4} & {\bf q^3} & {\bf q^2} & {\bf q } & {\bf 1}
\end{array}
\right]\]

C(2,8): 
\[ \left[
\begin{array}{cccccccccccc} r:         & 1      & 2   & 3   & 4   & 5 
&  6 & 7   & 8   & 9   & 10  & 11      \\
E_r:  & {\bf q^{12}} & {\bf q^{11}} &{\bf  q^{10}} & {\bf q^9} & {\bf 
q^8}
          & {\bf q^7} & {\bf q^6}  & q^{10} & q^9 & q^8 &q^7
\end{array}
\right]\]

\[\left[
\begin{array}{ccccccccccccccccc} 
12   & 13 & 14  & 15  & 16  & 17  & 18  & 19 &  20 & 21 & 22  & 23  & 
24  & 25  & 26  & 27  & 28 \\
  q^6  & q^5 & q^8 & q^4 & q^7 & q^6 & q^5 & q^4 &  q^3  & q^2 & {\bf
  q^6 }& {\bf q^5} & {\bf q^4} & {\bf q^3} & {\bf q^2} & {\bf q} &
  {\bf 1} 
\end{array}
\right]\]

%\begin{tiny}
C(3,6): 
\[ \left[
\begin{array}{ccccccccccc} r:         & 1      & 2   & 3   &
  4 & 5 & 6 & 7  & 8   & 9   & 10  \\
E_r:  & {\bf q^{9}} & {\bf q^{8}} & {\bf q^{7}} & {\bf q^6} & q^7 & 
q^5 & q^6 &  q^5 & q^4 &
q^3 
\end{array}
\right]\]

\[ \left[
\begin{array}{cccccccccc} 11 & 12   & 13 & 14  & 15  
& 16  & 17  & 18  & 19 & 20 \\
q^6   &  q^5  & q^4 & q^3 & q^4 & q^2 &{\bf q^3} & {\bf q^2}
&{\bf q } & {\bf 1 } 
\end{array}
\right]\]

%\end{tiny}
The expressions in boldface indicate values where $E_r=\Delta_r$ 
because of 
Theorem \ref{highweights} (to the right)
or Proposition \ref{Nogin} (to the left).
The expressions not in boldface contribute to upper bounds for ``the 
true values'' $d_r$, when adding monomials from left.

%\vfill\eject

As further illustration we give some more detailed tables of Schubert 
unions:

$G(2,5):$

\begin{tabular}{|c|c|c|c|c|c|} \hline
$$U$$ & Span & Krull & $M_U$ & number of points & Maximal \\ \hline
$\emptyset$ & $0$   & $-1$   & $\emptyset$ & $0$ & Yes \\ \hline
$(1,2)$ & $1$   & $0$ & $\{1\}$ & $1$ & Yes \\ \hline
$(1,3)$ & $2$   & $1$ & $\{2\}$ & $q+1$ & Yes \\ \hline
$(1,4)$ & $3$   & $2$ & $\{3\}$ & $q^2+q+1$ & Yes \\ \hline
$(1,5)$ & $4$   & $3$ & $\{4\}$ & $q^3+q^2+q+1$ & Yes \\ \hline
$(2,3)$ & $3$   & $2$ & $\{1,2\}$ & $q^2+q+1$ & Yes \\ \hline
$(1,4) \cup (2,3)$ & $4$   & $2$ & $\{1,3\}$ & $2q^2+q+1$ & No \\ 
\hline
$(1,5) \cup (2,3)$ & $5$   & $3$ & $\{1,4\}$ & $q^3+2q^2+q+1$ & Yes 
\\ \hline
$(2,4)$ & $5$   & $3$ & $\{2,3\}$ & $q^3+2q^2+q+1$ & Yes \\ \hline
$(1,5) \cup (2,4)$ & $6$   & $3$ & $\{2,4\}$ & $2q^3+2q^2+q+1$ & No \\
\hline
$(2,5)$ & $7$   & $4$ & $\{3,4\}$ & $q^4+2q^3+2q^2+q+1$ & Yes \\ 
\hline
$(3,4)$ & $6$   & $4$ & $\{1,2,3\}$ & $q^4+q^3+2q^2+q+1$ & Yes  \\ 
\hline
$(1,5) \cup (3,4)$ & $7$   & $4$ & $\{1,2,4\}$ & $q^4+2q^3+2q^2+q+1$ 
& Yes \\
\hline
$(2,5) \cup (3,4)$ & $8$   & $4$ & $\{1,3,4\}$ & $2q^4+2q^3+2q^2+q+1$ 
& Yes \\
\hline
$(3,5)$ & $9$   & $5$ & $\{2,3,4\}$ & $q^5+2q^4+2q^3+2q^2+q+1$ & Yes 
\\
\hline
$(4,5)$ & $10$   & $6$ & $\{1,2,3,4\}$ & $q^6+q^5+2q^4+2q^3+2q^2+q+1$
& Yes \\
\hline
\end{tabular}
\eject
In this table we have listed all non-trivial Schubert unions (for a 
fixed flag) 
for $G(2,5)$. Below we give similar tables for $G(2,6)$ and $G(3,6)$.
In the column to the right we indicate whether the Schubert union in
question has the maximum possible of points among the Schubert unions
of that spanning dimension. The (affine) spanning dimension is given
in the column marked ``Span''.

\vspace{.4cm}
%
%$G(2,4):$
%\begin{tabular}{|c|c|c|c|c|c|} \hline
%$$U$$ & Span & Krull & $M_U$ & number of points & Maximal\\ \hline
%$\emptyset$ & $-1$   & $0$   & $\emptyset$ & $0$ & Yes \\ \hline
%$(1,2)$ & $1$   & $0$ & $\{1\}$ & $1$ & Yes \\ \hline
%$(1,3)$ & $2$   & $1$ & $\{2\}$ & $q+1$ & Yes \\ \hline
%$(1,4)$ & $3$   & $2$ & $\{3\}$ & $q^2+q+1$ & Yes \\ \hline
%$(2,3)$ & $3$   & $2$ & $\{1,2\}$ & $q^2+q+1$ & Yes \\ \hline
%$(1,4) \cup (2,3)$ & $4$   & $2$ & $\{1,3\}$ & $2q^2+q+1$ & Yes \\
%\hline
%$(2,4)$ & $5$   & $3$ & $\{2,3\}$ & $q^3+2q^2+q+1$ & Yes \\ \hline
%$(3,4)$ & $6$   & $4$ & $\{1,2,3\}$ & $q^4+q^3+2q^2+q+1$ & Yes \\
%\hline
%\end{tabular}
%\bigskip

The dual of a Schubert
union with a given $M_U$ in $G(2,m)$ is the Schubert union $V$ with
$M_V=\{1,...,m-1\}-M_U$. 
%The same remark of course also applies to the
%following two tables for $G(2,5)$ and $G(2,6)$. 
We also remark that
even
in the cases where $m \choose 2$ is even, there can be no self-dual 
Schubert unions in $G(2,m)$ (of spanning dimension 
$m \choose 2$$/2$), since a set $M_U$ is never equal to its own
complement. We shall se below that the situation may be different
for $G(l,m)$ with $l=3$.

\bigskip

%\eject
$G(2,6):$

%\begin{tiny}
\begin{tabular}{|c|c|c|c|c|c|} \hline
$$U$$ & Span & Krull & $M_U$ & number of points & Max.\\ \hline
$\emptyset$ & $0$   & $-1$   & $\emptyset$ & $0$ & Yes \\ \hline
$(1,2)$ & $1$   & $0$ & $\{1\}$ & $1$  & Yes \\ \hline
$(1,3)$ & $2$   & $1$ & $\{2\}$ & $q+1$  & Yes \\ \hline
$(1,4)$ & $3$   & $2$ & $\{3\}$ & $q^2+q+1$ & Yes  \\ \hline
$(1,5)$ & $4$   & $3$ & $\{4\}$ & $q^3+q^2+q+1$  & Yes \\ \hline
$(1,6)$ & $5$   & $4$ & $\{5\}$ & $q^4+q^3+q^2+q+1$ & Yes  \\ \hline
$(2,3)$ & $3$   & $2$ & $\{1,2\}$ & $q^2+q+1$ & Yes \\ \hline
$(1,4) \cup (2,3)$ & $4$   & $2$ & $\{1,3\}$ & $2q^2+q+1$ &
No  \\ \hline
$(1,5) \cup (2,3)$ & $5$   & $3$ & $\{1,4\}$ & $q^3+2q^2+q+1$ &
No  \\ \hline
$(1,6) \cup (2,3)$ & $6$   & $4$ & $\{1,5\}$ & $q^4+q^3+2q^2+q+1$ & 
Yes \\
\hline
$(2,4)$ & $5$   & $3$ & $\{2,3\}$ & $q^3+2q^2+q+1$ & No \\ \hline
$(1,5) \cup (2,4)$ & $6$   & $3$ & $\{2,4\}$ & $2q^3+2q^2+q+1$ &
No \\
\hline
$(1,6) \cup (2,4)$ & $7$   & $4$ & $\{2,5\}$ & $q^4+2q^3+2q^2+q+1$ & 
Yes \\
\hline
$(2,5)$ & $7$   & $4$ & $\{3,4\}$ & $q^4+2q^3+2q^2+q+1$ & Yes \\ 
\hline
$(1,6) \cup (2,5)$ & $8$   & $4$ & $\{3,5\}$ & $2q^4+2q^3+2q^2+q+1$  &
Yes \\
\hline
$(2,6)$ & $9$   & $5$ & $\{4,5\}$ & $q^5+2q^4+2q^3+2q^2+q+1$ & Yes \\
\hline
$(3,4)$ & $6$   & $4$ & $\{1,2,3\}$ & $q^4+q^3+2q^2+q+1$ & Yes \\ 
\hline
$(1,5) \cup (3,4)$ & $7$   & $4$ & $\{1,2,4\}$ & $q^4+2q^3+2q^2+q+1$ &
Yes \\
\hline
$(1,6) \cup (3,4)$ & $8$   & $4$ & $\{1,2,5\}$ & $2q^4+2q^3+2q^2+q+1$ 
& Yes \\
\hline
$(2,5) \cup (3,4)$ & $8$   & $4$ & $\{1,3,4\}$ & 
$2q^4+2q^3+2q^2+q+1$& Yes  \\
\hline
$(1,6) \cup  (2,5) \cup (3,4)$ & $9$   & $4$ & $\{1,3,5\}$ &
$3q^4+2q^3+2q^2+q+1$  & No \\
\hline
$(2,6) \cup (3,4)$ & $10$   & $5$ & $\{1,4,5\}$ & 
$q^5+3q^4+2q^3+2q^2+q+1$ &
No \\
\hline
$(3,5)$ & $9$   & $5$ & $\{2,3,4\}$ & $q^5+2q^4+2q^3+2q^2+q+1$ & Yes
\\
\hline
$(1,6) \cup  (3,5)$ & $10$   & $5$ & $\{2,3,5\}$ & 
$q^5+3q^4+2q^3+2q^2+q+1$ &
No \\
\hline
$(2,6) \cup  (3,5)$ & $11$   & $5$ & $\{2,4,5\}$ & 
$2q^5+3q^4+2q^3+2q^2+q+1$ &
No  \\
\hline
$(3,6)$ & $12$   & $6$ & $\{3,4,5\}$ & $q^6+2q^5+3q^4+2q^3+2q^2+q+1$ &
Yes \\
\hline
$(4,5)$ & $10$   & $6$ & $\{1,2,3,4\}$ & $q^6+q^5+2q^4+2q^3+2q^2+q+1$
& Yes \\
\hline
$(1,6) \cup  (4,5)$ & $11$   & $6$ & $\{1,2,3,5\}$ & 
$q^6+q^5+3q^4+2q^3+2q^2+q+1$ & Yes \\
\hline
$(2,6) \cup  (4,5)$ & $12$   & $6$ & $\{1,2,4,5\}$ &
$q^6+2q^5+3q^4+2q^3+2q^2+q+1$  & Yes \\
\hline
$(3,6) \cup  (4,5)$ & $13$   & $6$ & $\{1,3,4,5\}$ &
$2q^6+2q^5+3q^4+2q^3+2q^2+q+1$  & Yes \\
\hline
$(3,6) \cup  (4,5)$ & $14$   & $7$ & $\{2,3,4,5\}$ &
$q^7+2q^6+2q^5+3q^4+2q^3+2q^2+q+1$  & Yes \\
\hline
$(5,6)$ & $15$   & $8$ & $\{1,2,3,4,5\}$ & $n$ & Yes  \\
\hline
\end{tabular}
%\end{tiny}
\bigskip
\eject
%The corresponding (smaller) table for $G(2,4)$ and $G(2,5)$, given
%earlier, can be 
%derived from the table for $G(2,6)$, roughly speaking by only
%focusing on those rows
%where $M_U$ is a subset of $\{1,2,3\}$, or $\{1,2,3,4\}$,
%respectively. A little caution is necessary, though. For $G(2,4)$
%one quickly sees that all $8$ Schubert unions are maximal for their
%spanning dimensions. For $G(2,5)$ all $16$ unions, except those with 
%$M_U=\{1,3\}$ or $\{2,4\}$ are maximal. For $G(2,6)$ we for example see that
%unions $U$ with $M_U=\{2,3\}, \{1,4\}, \{2,4\}$ are not maximal.
%So it is not an ``intrinsic'' property of a Schubert union whether it
%is maximal for its spanning dimension. It depends on the Grassmann
%variety, in which it sits.

$G(3,6):$

\begin{tabular}{|c|c|c|c|c|c|c|} \hline
$$U$$ & Span & Krull  & Max.& Number of points \\ \hline
$\emptyset$ & $0$   & $-1$  & Yes & $0$\\ \hline
$(1,2,3)$ & $1$   & $0$  & Yes & $1$ \\ \hline
$(1,2,4)$ & $2$   & $1$  & Yes & $q+1$ \\ \hline
$(1,2,5)$ & $3$   & $2$ & Yes & $q^2+q+1$ \\ \hline
$(1,2,6)$ & $4$   & $3$ & Yes & $q^3+q^2+q+1$ \\ \hline
$(1,3,4)$ & $3$   & $2$ & Yes &  $q^2+q+1$\\ \hline
$(1,3,5)$ & $5$   & $3$ & Yes &  $q^3+2q^2+q+1$\\ \hline
$(1,3,6)$ & $7$   & $4$ & Yes &  $q^4+2q^3+2q^2+q+1$  \\ \hline
$(1,4,5)$ & $6$   & $4$  & Yes &  $q^4+q^3+2q^2+q+1$ \\ \hline
$(1,4,6)$ & $9$   & $5$  & Yes &  $q^5+2q^4+2q^3+2q^2+q+1$ \\ \hline
$(1,5,6)$ & $10$   & $6$ & Yes &  $q^6+q^5+2q^4+2q^3+2q^2+q+1$  \\ 
\hline
$(2,3,4)$ & $4$   & $3$ & Yes & $q^3+q^2+q+1$ \\ \hline
$(2,3,5)$ & $7$   & $4$ & Yes &  $q^4+2q^3+2q^2+q+1$  \\ \hline
$(2,3,6)$ & $10$   & $5$ & No   &  $q^5+2q^4+3q^3+2q^2+q+1$ \\ \hline
$(2,4,5)$ & $9$   & $5$ & Yes    &  $q^5+2q^4+2q^3+2q^2+q+1$ \\ \hline
$(2,4,6)$ & $14$   & $6$  & No   &  $q^6+3q^5+3q^4+3q^3+2q^2+q+1$ \\ 
\hline
$(2,5,6)$ & $16$   & $7$ & Yes  &  $q^7+2q^6+3q^5+3q^4+3q^3+2q^2+q+1$ 
\\ \hline
$(3,4,5)$ & $10$   & $6$ & Yes  &  $q^6+q^5+2q^4+2q^3+2q^2+q+1$\\ 
\hline
$(3,4,6)$ & $16$   & $7$ & Yes   &  
$q^7+2q^6+3q^5+3q^4+3q^3+2q^2+q+1$\\ \hline
$(3,5,6)$ & $19$   & $8$ & Yes   &  
$q^8+2q^7+3q^6+3q^5+3q^4+3q^3+2q^2+q+1$\\ \hline
$(4,5,6)$ & $20$   & $9$ & Yes   &  
$q^9+q^8+2q^7+3q^6+3q^5+3q^4+3q^3+2q^2+q+1$ \\ \hline
$(1,2,5) \cup (1,3,4)$ & $4$   & $2$ & No & $2q^2+q+1$ \\ \hline
$(1,2,5) \cup (2,3,4)$ & $5$   & $3$ & Yes & $q^3+2q^2+q+1$ \\ \hline
$(1,2,6) \cup (1,3,4)$ & $5$   & $3$ & Yes & $q^3+2q^2+q+1$ \\ \hline
$(1,2,6) \cup (1,3,5)$ & $6$   & $3$ & No  & $2q^3+2q^2+q+1$ \\ \hline
$(1,2,6) \cup (1,4,5)$ & $7$   & $4$ & Yes & $q^4+2q^3+2q^2+q+1$\\ 
\hline
$(1,2,6) \cup (2,3,4)$ & $6$   & $3$ & No & $2q^3+2q^2+q+1$ \\ \hline
$(1,2,6) \cup (2,3,5)$ & $8$   & $4$ & No  & $q^4+3q^3+2q^2+q+1$\\ 
\hline
$(1,2,6) \cup (2,4,5)$ & $10$   & $5$ & No  & 
$q^5+2q^4+3q^3+2q^2+q+1$\\ \hline
$(1,3,5) \cup (2,3,4)$ & $6$   & $3$ & No & $2q^3+2q^2+q+1$ \\ \hline
$(1,3,6) \cup (1,4,5)$ & $8$   & $4$ & Yes  & $2q^4+2q^3+2q^2+q+1$ \\ 
\hline
$(1,3,6) \cup (3,4,5)$ & $12$   & $6$ & Yes & 
$q^6+q^5+3q^4+3q^3+2q^2+q+1$  \\ \hline
$(1,3,6) \cup (2,3,4)$ & $8$   & $4$ & No   & $q^4+3q^3+2q^2+q+1$\\ 
\hline
$(1,3,6) \cup (2,3,5)$ & $9$   & $4$ & No & $2q^4+3q^3+2q^2+q+1$ \\ 
\hline
$(1,3,6) \cup (2,4,5)$ & $11$   & $5$ & No & $q^5+3q^4+3q^3+2q^2+q+1$ 
\\ \hline
$(1,4,6) \cup (2,3,4)$ & $10$   & $5$ & No & 
$q^5+2q^4+3q^3+2q^2+q+1$\\ \hline
$(1,4,6) \cup (2,3,5)$ & $11$   & $5$ & No & $q^5+3q^4+3q^3+2q^2+q+1$ 
\\ \hline
$(1,4,6) \cup (2,3,6)$ & $12$   & $5$ & No  & 
$2q^5+3q^4+3q^3+2q^2+q+1$\\ \hline
$(1,4,6) \cup (2,4,5)$ & $12$   & $5$ & No & 
$2q^5+3q^4+3q^3+2q^2+q+1$ \\ \hline
$(1,4,6) \cup (3,4,5)$ & $13$   & $6$ & Yes & 
$q^6+2q^5+3q^4+3q^3+2q^2+q+1$ \\ \hline
$(1,5,6) \cup (2,3,4)$ & $11$   & $6$ & Yes &  
$q^6+q^5+2q^4+3q^3+2q^2+q+1$  \\ \hline
$(1,5,6) \cup (2,3,5)$ & $12$   & $6$ & Yes & 
$q^6+q^5+3q^4+3q^3+2q^2+q+1$ \\ \hline

\end{tabular}

\begin{tabular}{|c|c|c|c|c|c|c|} \hline
$$U$$ & Span & Krull  & Max.& Number of points\\ \hline
$(1,5,6) \cup (2,3,6)$ & $13$   & $6$ & Yes &
$q^6+2q^5+3q^4+3q^3+2q^2+q+1$ \\ \hline
$(1,5,6) \cup (2,4,5)$ & $13$   & $6$ & Yes & 
$q^6+2q^5+3q^4+3q^3+2q^2+q+1$ \\ \hline
$(1,5,6) \cup (2,4,6)$ & $15$   & $6$ & Yes & 
$2q^6+3q^5+3q^4+3q^3+2q^2+q+1$  \\ \hline
$(1,5,6) \cup (3,4,5)$ & $14$   & $6$ & Yes  & 
$2q^6+2q^5+3q^4+3q^3+2q^2+q+1$ \\ \hline
$(1,5,6) \cup (3,4,6)$ & $17$   & $7$ & Yes & 
$q^7+3q^6+3q^5+3q^4+3q^3+2q^2+q+1$ \\ \hline
$(2,3,6) \cup (2,4,5)$ & $12$   & $5$ & No  & 
$2q^5+3q^4+3q^3+2q^2+q+1$\\ \hline
$(2,3,6) \cup (3,4,5)$ & $13$   & $6$ & Yes & 
$q^6+2q^5+3q^4+3q^3+2q^2+q+1$ \\ \hline
$(2,4,6) \cup (3,4,5)$ & $15$   & $6$ & Yes & 
$2q^6+3q^5+3q^4+3q^3+2q^2+q+1$ \\ \hline
$(2,5,6) \cup (3,4,5)$ & $17$   & $7$ & Yes  & 
$q^7+3q^6+3q^5+3q^4+3q^3+2q^2+q+1$ \\ \hline
$(2,5,6) \cup (3,4,6)$ & $18$   & $7$ & Yes  & 
$2q^7+3q^6+3q^5+3q^4+3q^3+2q^2+q+1$ \\ \hline
$(1,4,5) \cup (2,3,5)$ & $8$   & $4$ & Yes & $2q^4+2q^3+2q^2+q+1$ \\ 
\hline
$(1,4,5) \cup (2.3,4)$ & $7$   & $4$ & Yes & $q^4+2q^3+2q^2+q+1$  \\ 
\hline
$(1,2,6) \cup (3,4,5)$ & $11$   & $6$ & Yes & 
$q^6+q^5+2q^4+3q^3+2q^2+q+1$ \\ \hline
$(1,4,5) \cup (2,3,6)$ & $11$   & $5$ & No & 
$q^5+3q^4+3q^3+2q^2+q+1$  \\ \hline
$(1,3,6) \cup (1,4,5) \cup (2,3,4)$ & $9$   & $4$ & No & 
$2q^4+3q^3+2q^2+q+1$ \\ \hline
$(1,4,6) \cup (2,3,6) \cup (2,4,5)$ & $13$   & $5$ & No & 
$3q^5+3q^4+3q^3+2q^2+q+1$  \\ \hline
$(1,4,6) \cup (2,3,6) \cup (3,4,5)$ & $14$   & $6$ & No & 
$q^6+3q^5+3q^4+3q^3+2q^2+q+1$  \\ \hline
$(1,5,6) \cup (2,3,6) \cup (2,4,5)$ & $14$   & $6$ & No & 
$q^6+3q^5+3q^4+3q^3+2q^2+q+1$  \\ \hline
$(1,5,6) \cup (2,3,6) \cup (3,4,5)$ & $15$   & $6$ & Yes  & 
$2q^6+3q^5+3q^4+3q^3+2q^2+q+1$ \\ \hline
$(1,5,6) \cup (2,4,6) \cup (3,4,5)$ & $16$   & $6$ & No  & 
$3q^6+3q^5+3q^4+3q^3+2q^2+q+1$ \\ \hline
$(1,2,6) \cup (1,3,5) \cup (2,3,4)$ & $7$   & $3$ & No  & 
$3q^3+2q^2+q+1$ \\ \hline
$(1,2,6) \cup (1,4,5) \cup (2,3,5)$ & $9$   & $4$ & No & 
$2q^4+3q^3+2q^2+q+1$ \\ \hline
$(1,2,6) \cup (1,4,5) \cup (2,3,4)$ & $8$   & $4$ & No & 
$q^4+3q^3+2q^2+q+1$ \\ \hline
$(1,3,6) \cup (1,4,5) \cup (2,3,5)$ & $10$   & $4$ & No & 
$3q^4+3q^3+2q^2+q+1$ \\ \hline

\end{tabular}

\bigskip

All the tables made so far have been produced, mainly by using 
Corollary \ref{int-un}. The Krull dimension 
(of the component cycle with biggest such dimension) is of course 
equal to
the degree of $g_U(q)$, interpreted as a polynomial in $q$,
and this is the polynomial appearing in the column marked ``number of 
points''.
Moreover it is well known that the Krull dimension of a Schubert cycle
$S_{(a_1,...,a_l)}$ is $a_1+a_2+...+a_l-\frac{l(l+1)}{2}$, so the 
Krull dimension can be ``read off `` both from the leftmost and the
rightmost column.

In the table for $G(3,6)$ above study the $16$ rows corresponding to 
unions of
cycles $S_{(a,b,c)}$ with $c \le 5$. This gives rise to the 
corresponding
table for $G(3,5)$. But this is isomorphic to $G(2,5)$. It is an
amusing exercise to translate all unions in  $G(3,5)$ to corresponding
ones in $G(2,5)$ and check that the relevant columns of the tables
coincide. 

\begin{rem} \label{order}
{\rm Given two Schubert unions $U_1, U_2$ with corresponding
  polynomials
$g_{U_1}(q)$ and  $g_{U_2}(q)$. The issue of which of the two that
  gives the highest value for given $q$ is in principle a different
  one, for each $q$. On the other hand, if we order the Schubert 
unions 
lexicographically with 
respect to $g_U$, then as remarked in Section \ref{algo}, it is clear 
that 
the lexicographic order is the
same as the ``number of point''-order for all large enough $q$.
In all the examples we have seen up to now, it is clear by inspection
  that these orders are the same for all prime powers $q$.
Hence the ``Yes'' and ``No'' in the ``Max.'' column can be interpreted
in two ways simultaneously (counting points, and ordering with 
respect to $g_U$).
}
\end{rem}
\eject
Table of dual pairs of Schubert unions for $G(3,6)$:
\bigskip

\begin{tabular}{|c|c|c|c|c|} \hline
$$U$$ & Span & Dual Schubert union  & Max. \\ \hline
$\emptyset$ & $0$   & $(4,5,6)$  & Yes \\ \hline
$(1,2,3)$ & $1$   & $(3,5,6)$  & Yes \\ \hline
$(1,2,4)$ & $2$   & $(2,5,6) \cap (3,4,6)$  & Yes  \\ \hline
$(1,2,5)$ & $3$   & $(1,5,6) \cap (3,4,6)$ & Yes  \\ \hline
$(1,3,4)$ & $3$   & $(2,5,6) \cap (3,4,5)$ & Yes \\ \hline
$(1,2,6)$ & $4$   & $(3,4,6)$ & Yes  \\ \hline
$(2,3,4)$ & $4$   & $(2,5,6)$ & Yes  \\ \hline
$(1,2,5) \cup (1,3,4)$ & $4$   & $(1,5,6)\cup (2,4,6) \cup (3,4,5)$ & 
No \\ \hline
$(1,3,5)$ & $5$   & $(1,5,6)\cup (2,3,6) \cup (3,4,5)$ & Yes \\ \hline
$(1,2,5) \cup (2,3,4)$ & $5$   & $(1,5,6)\cup (2,4,6)$ &
Yes \\ \hline
$(1,2,6) \cup (1,3,4)$ & $5$   & $(2,4,6)\cup (3,4,5)$ &
Yes \\ \hline
$(1,4,5)$ & $6$   & $(1,5,6) \cup ((3,4,5)$  & Yes  \\ \hline
$(1,2,6) \cup (1,3,5)$ & $6$   & $(1,4,6)\cup (2,3,6) \cup (3,4,5)$ &
No \\ \hline
$(1,2,6) \cup (2,3,4)$ & $6$   & $(2,4,6)$ &
No \\ \hline
$(1,3,5) \cup (2,3,4)$ & $6$   & $(1,5,6)\cup (2,3,6) \cup ((2,4,5)$ &
No \\ \hline
$(1,3,6)$ & $7$   & $(2,3,6) \cup (3,4,5)$ & Yes  \\ \hline
$(2,3,5)$ & $7$   & $(1,5,6) \cup (2,3,6)$ & Yes   \\ \hline
$(1,2,6) \cup (1,4,5)$ & $7$   & $(1,4,6)\cup (3,4,5) $ &
Yes \\ \hline
$(1,4,5) \cup (2,3,4)$ & $7$   & $(1,5,6)\cup (2,4,5)$ &
Yes \\ \hline
$(1,2,6) \cup (1,3,5) \cup (2,3,4)$ & $7$   & $(1,4,6)\cup (2,3,6) 
\cup (2,4,5)$ &
No \\ \hline
$(1,3,6) \cup (1,4,5)$ & $8$   & $(1,3,6)\cup (3,4,5)$ &
Yes \\ \hline
$(1,4,5) \cup (2,3,5)$ & $8$   & $(1,5,6)\cup (2,3,5)$ &
Yes \\ \hline
$(1,2,6) \cup (2,3,5)$ & $8$   & $(1,4,6)\cup (2,3,6)$ &
No \\ \hline
$(1,3,6) \cup (2,3,4)$ & $8$   & $(2,3,6)\cup (2,4,5)$ &
No \\ \hline
$(1,2,6) \cup (1,4,5) \cup (2,3,4)$ & $8$   & $(1,4,6)\cup (2,4,5)$ &
No \\ \hline
$(1,4,6)$ & $9$   & $(1,2,6) \cup (3,4,5)$  & Yes  \\ \hline
$(2,4,5)$ & $9$   & $(1,5,6) \cup (2,3,4)$ & Yes    \\ \hline
$(1,3,6) \cup (2,3,5)$ & $9$   & $(1,4,5)\cup (2,3,6)$ &
No \\ \hline
$(1,3,6) \cup (1,4,5) \cup (2,3,4)$ & $9$   & $(1,3,6)\cup (2,4,5)$ &
No \\ \hline
$(1,2,6) \cup (1,4,5) \cup (2,3,5)$ & $9$   & $(1,4,6)\cup (2,3,5)$ &
No \\ \hline
$(1,5,6)$ & $10$   & $(3,4,5)$ & Yes  \\ \hline
$(2,3,6)$ & $10$   & $(2,3,6)$ & No    \\ \hline
%$(3,4,5)$ & $10$   & $(1,5,6)$ & Yes  \\ \hline
$(1,2,6) \cup (2,4,5)$ & $10$   & $(1,4,6)\cup (2,3,4)$ &
No \\ \hline
$(1,3,6) \cup (1,4,5) \cup (2,3,5)$ & $10$   & $(1,3,6)\cup (1,4,5)
\cup (2,3,5)$ &
No \\ \hline
\end{tabular}

\bigskip
Above we also present a table of
Schubert unions and their associated dual Schubert unions for
$G(3,6)$.
We have used the methods described in Subsection \ref{techn} to make
the table.
All Schubert unions with spanning dimension at most $9$ can be found
in the left half of the table, and  unions with spanning dimension at 
least $11$ can be found on
the right side (as duals). For spanning dimension $10$  all
$6$ unions are listed on at least one side.
\begin{rem} \label{3rem}
 {\rm 
(i) The table  reveals
a situation 
different from the case $l=2$ and $m \choose l$ even, where no
Schubert union is self-dual.
Here we see that both $(2,3,6)$ and $(1,3,6) \cup (1,4,5) \cup
(2,3,5)$ are self-dual Schubert unions. 

(ii) Another fact that can be found from the table is that the 
conclusion of
Corollary \ref{moredualsigma} fails for $l=3$. We see that the dual of
$S_{(1,3,5)}$ is the proper triple union 
$S_{(1,5,6)} \cup S_{(2,3,6)} \cup S_{(3,4,5)}$ (and vice versa).
%We encourage the interested reader to reconstruct this situation,
%and the selfduality described in (i),  
%by playing with cubes. 

(iii) For this Grassmann varieties we have described in the tables 
above, a
Schubert union
has a maximal number of points, given its spanning dimension, if and 
only if
its dual union enjoys the same property, so the ``Yes'' and ``No'' in
the ``Max.''-column of the last table apply to the left and right half
of  the table simultaneously. The same nice property also holds for  
$(l,m)=(2,7)$ 
and $(2,8)$, but for reasons of space we do not give the full tables
here, from which the shorter lists of the $E_r$ at the start of this 
section were deduced.}

%This remark implies that the answer to question (Q8) of
%Section \ref{w67} is affirmative.}

\end{rem}


\begin{thebibliography}{[E-L-M-S]}

\bibitem[GL]{GL}  Ghorpade, S, Lachaud, G. \textit{ Higher weights
of Grassmann Codes}, in \textit{Coding Theory, Cryptography, and 
Related
Areas} (Guanajuoto, 1998), Springer Verlag, Berlin/Heidelberg, 
 122-31 (2000).

\bibitem[Fu]{Fu}  W. Fulton,  \textit{Young Tableaux}, 
 Student Texts, 35, London Math. Soc. (1991).

\bibitem[FH]{FH}  W. Fulton, J. Harris,   \textit{Representation Theory}, 
Graduate Texts in Mathematics, 129, Springer Verlag (1991).

\bibitem[GT]{GT} Ghorpade, S., Tsfasman, M. \textit{Schubert 
Varieties,
Linear Codes, and Enumerative Combinatorics}, Preprint (2003).

\bibitem[HC]{HC}  H. Chen,  \textit{ On the minimum distance of 
Schubert 
codes}, IEEE Trans. of Inform. Theory, 46, 1535-38 (2000).

\bibitem[H]{H}  R. Hartshorne,  \textit{Algebraic Geometry}, 
Graduate Texts in Mathematics, Springer Verlag, (1977).

\bibitem[HT]{HT}  J.P. Hirschfeld, J. A.Thas,  \textit{General Galois 
Geometries}, Oxford: Clarendon Press(1991).

\bibitem[MV]{MV}  Elisa Montanucci, Rita Vinscenti
  \textit{Characterization of linear codes and classical varieties}, 
Preprint (2004).

\bibitem[N]{N} Nogin, D.Yu., \textit{Codes associated to 
Grassmannians}, 
in \textit{Arithmetic Geometry and Coding theory} (Luminy 1993),
R. Pellikaan, M. Perret, S.G. Vladut, Eds. Walter de Gruyter,
Berlin/New York,  145-54 (1996).

\bibitem[Sp]{Sp}  T. A. Springer,
  \textit{Linear Algebraic Groups}, Progress in Mathematics, 9, Birkh\"auser
 (1981).

\end{thebibliography}
\end{document}